\documentclass[%
  a4paper,
  onecolumn,
]{mypreprint}

\usepackage{amsmath}
\usepackage{amssymb}
\usepackage{amsthm}

\usepackage[english]{babel}

\usepackage{tikz}
  \usetikzlibrary{positioning}

\usepackage{pgfplots}
  \pgfplotsset{compat = 1.13}
  \usepgfplotslibrary{external}
  \tikzset{external/system call = {%
    pdflatex \tikzexternalcheckshellescape
      -halt-on-error
      -interaction=batchmode
      -jobname "\image" "\texsource"}}
  \tikzexternaldisable

\usepackage{filemod}

\usepackage[linesnumbered, lined, algoruled]{algorithm2e}

\crefname{algocf}{alg.}{algs.}
\Crefname{algocf}{Algorithm}{Algorithms}
\Crefname{problem}{Problem}{Problems}

\usepackage{enumitem}

\usepackage{todonotes}
\usepackage{cleveref}

\newcommand{\Hinf}{\ensuremath{\mathcal{H}_{\infty}}}
\newcommand{\trans}{\ensuremath{\mkern-1.5mu\mathsf{T}}}

\newcommand{\Xinf}{\ensuremath{X_{\infty}}}

\newcommand{\hA}{\ensuremath{\widehat{A}}}
\newcommand{\hB}{\ensuremath{\widehat{B}}}
\newcommand{\hC}{\ensuremath{\widehat{C}}}
\newcommand{\hE}{\ensuremath{\widehat{E}}}
\newcommand{\hX}{\ensuremath{\widehat{X}}}

\newcommand{\cP}{\ensuremath{\mathcal{P}}}

\newcommand{\R}{\ensuremath{\mathbb{R}}}
\newcommand{\C}{\ensuremath{\mathbb{C}}}

\newcommand{\cR}{\ensuremath{\mathcal{R}}}
\newcommand{\tA}{\ensuremath{\skew4\widetilde{A}}}
\newcommand{\tE}{\ensuremath{\skew4\widetilde{E}}}
\newcommand{\tB}{\ensuremath{\skew4\widetilde{B}}}
\newcommand{\tC}{\ensuremath{\skew4\widetilde{C}}}
\newcommand{\tn}{\ensuremath{\tilde{n}}}

\newcommand{\mnum}[3]{\ensuremath{#1\texttt{e#2}#3}}

\newcommand\Tstrut{\rule{0pt}{3.4ex}}
\newcommand\Bstrut{\rule[-1.3ex]{0pt}{0pt}}

\theoremstyle{plain}\newtheorem{proposition}{Proposition}
\theoremstyle{definition}\newtheorem{problem}{Problem}
\theoremstyle{definition}\newtheorem{remark}{Remark}
\theoremstyle{definition}

\definecolor{matlabBlue}{HTML}{0072BD}
\definecolor{matlabOrange}{HTML}{D95319}
\definecolor{matlabYellow}{HTML}{EDB120}
\definecolor{matlabPurple}{HTML}{7E2F8E}
\definecolor{matlabGreen}{HTML}{77AC30}
\definecolor{matlabLightBlue}{HTML}{4DBEEE}
\definecolor{matlabRed}{HTML}{A2142F}

\tikzstyle{ri} = [
  matlabBlue,
  line width   = 1.25pt,
  mark         = *,
  mark size    = 2pt,
  mark options = {solid}
]
\tikzstyle{inner1} = [
  matlabOrange,
  line width   = 1.25pt,
  mark         = square*,
  mark size    = 2pt,
]
\tikzstyle{inner2} = [
  matlabPurple,
  line width   = 1.25pt,
  mark         = triangle*,
  mark size    = 2pt,
]
\tikzstyle{inner3} = [
  matlabGreen,
  line width   = 1.25pt,
  mark         = diamond*,
  mark size    = 2pt,
]
\tikzstyle{tol} = [
  black,
  line width = 1.25pt,
  dashed
]

\begin{document}


\title{A low-rank solution method for Riccati equations with indefinite
  quadratic terms}

\author[$\ast$]{Peter Benner}
\affil[$\ast$]{Max Planck Institute for Dynamics of Complex Technical
  Systems, Sandtorstra{\ss}e 1, 39106 Magdeburg, Germany.
  \email{benner@mpi-magdeburg.mpg.de}, \orcid{0000-0003-3362-4103}
  \authorcr \itshape
  Otto von Guericke University Magdeburg, Faculty of Mathematics,
  Universit{\"a}tsplatz 2, 39106 Magdeburg, Germany.
  \email{peter.benner@ovgu.de}}
  
\author[$\dagger$]{Jan Heiland}
\affil[$\dagger$]{Max Planck Institute for Dynamics of Complex Technical
  Systems, Sandtorstra{\ss}e 1, 39106 Magdeburg, Germany.
  \email{heiland@mpi-magdeburg.mpg.de}, \orcid{0000-0003-0228-8522}
  \authorcr \itshape
  Otto von Guericke University Magdeburg, Faculty of Mathematics,
  Universit{\"a}tsplatz 2, 39106 Magdeburg, Germany.
  \email{jan.heiland@ovgu.de}}
  
\author[$\ddag$]{Steffen W. R. Werner}
\affil[$\ddag$]{Courant Institute of Mathematical Sciences, New York University,
  New York, NY 10012, USA.
  \email{steffen.werner@nyu.edu}, \orcid{0000-0003-1667-4862}}

\shorttitle{Low-rank Riccati iteration method}
\shortauthor{P. Benner, J. Heiland, S.~W.~R. Werner}
\shortdate{2021-11-12}
\shortinstitute{}
  

  
\abstract{%
  Algebraic Riccati equations with indefinite quadratic terms play an important
  role in applications related to robust controller design.
  While there are many established approaches to solve these in case of
  small-scale dense coefficients, there is no approach available to compute
  solutions in the large-scale sparse setting.
  In this paper, we develop an iterative method to compute low-rank
  approximations of stabilizing solutions of large-scale sparse continuous-time
  algebraic Riccati equations with indefinite quadratic terms.
  We test the developed approach for dense examples in comparison to other
  established matrix equation solvers, and investigate the applicability and
  performance in large-scale sparse examples.
}

\novelty{%
  We propose an iterative algorithm for computing low-rank approximations
  of stabilizing solutions of continuous-time algebraic Riccati equations with
  indefinite quadratic term.
  This is the first such approach that applies for general setups with
  large-scale sparse coefficient matrices.
}

\maketitle


\section{Introduction}%
\label{sec:intro}

Many concepts in systems and control theory are connected to solutions of
algebraic Riccati equations.
Prominent examples are the linear-quadratic regulator (LQR) and linear-quadratic
Gaussian (LQG) controller design~\cite{Son98, AndM90, Loc01, LanR95} and
corresponding model order reduction methods~\cite{morJonS83} as well as the
characterization of passivity and contractivity of input-output systems and
properties-preserving model reduction methods for
these~\cite{morDesP82, morOpdJ88, AndV72}.
They also appear, for example, in applications with differential
games~\cite{Del07, BasM17}.
In this paper, we consider algebraic Riccati equations with indefinite quadratic
terms of the form
\begin{align} \label{eqn:indric}
  A^{\trans} X E + E^{\trans} X A + E^{\trans} X \left( B_{1} B_{1}^{\trans} -
  B_{2} B_{2}^{\trans} \right) X E + C^{\trans} C & = 0,
\end{align}
with $A, E \in \R^{n \times n}$, $B_{1} \in \R^{n \times m_{1}}$,
$B_{2} \in \R^{n \times m_{2}}$, $C \in \R^{p \times n}$ and $E$ invertible,
and we are interested in symmetric positive semi-definite stabilizing solutions
of~\cref{eqn:indric}.
Equations of the form~\cref{eqn:indric} usually occur in $\Hinf$-control theory
and related topics, e.g., robust controller
design~\cite{McFG90, DoyGKF89, morMusG91}.

In the case of coefficients of small size $n$ and classical LQR/LQG Riccati
equations, i.e., $B_{1} = 0$ in~\cref{eqn:indric} such that the quadratic term
is negative semi-definite, there is a variety of different numerical approaches
to compute the stabilizing solution: Direct approaches that
compute eigenvalue decompositions of underlying Hamiltonian or even matrix
pencils~\cite{Lau79, AmmBM93, ArnL84}, iterative methods that work on the
underlying spectrum of the Hamiltonian matrix~\cite{morRob80, Ben97,
BenEQetal14}, or iterative approaches, which compute a sequence of matrices
converging to the stabilizing solution~\cite{Kle68, San74}.
The problem becomes more complicated in the opposite case with 
symmetric positive semi-definite quadratic term, i.e., in~\cref{eqn:indric}
setting $B_{2} = 0$, which occurs in bounded-real and positive-real
problems~\cite{morDesP82, morOpdJ88, AndV72}.
This already reduces the number of applicable methods.
While methods from the classical setting that do not take the definiteness
of the quadratic term into account can still be applied
here~\cite{Lau79, ArnL84, morRob80, Ben97}, only few, iterative Newton-type
methods, which converge to the desired solution, have been developed for this
problem class~\cite{BenS14, Var95b, Ben97a}.
The amount of applicable methods narrows down even further when considering
the general case of indefinite quadratic terms like in~\cref{eqn:indric}.
A new type of iterative method for the solution of~\cref{eqn:indric} was
developed in~\cite{LanFA07,LanFAetal08} to overcome accuracy problems of
classical solution approaches and the general lack of iterative methods for this
problem type.
This new method computes a sequence of matrices converging to the stabilizing
solution of~\cref{eqn:indric} by solving classical Riccati equations with
symmetric negative semi-definite quadratic terms.

However, in this paper, we will mainly focus on the case of large-scale sparse
coefficient matrices in~\cref{eqn:indric} arising, for example, from the
discretization of partial differential equations, with $n \gtrapprox 10^{5}$
and low-rank quadratic term and right-hand side such that
$m_{1}, m_{2}, p \ll n$.
This leads to another multitude of problems.
First of all, the dense computation methods from above cannot be applied anymore
since they may transform the original data and compute full solutions, which
easily becomes unfeasible in terms of memory requirements.
Also, these methods usually use matrix operations that cannot be efficiently
used in the large-scale sparse setting, which heavily increases the running
time of algorithms.
As in the dense case, there exists a variety of methods for the classical
LQR/LQG case ($B_{1} = 0$), for example, the low-rank Newton
method~\cite{BenLP08, Wei16} with an underlying low-rank alternating direction
implicit method~\cite{LiW02, BenLP08, BenKS13a, BenKS14b, Kue16} for solving the
occurring large-scale sparse Lyapunov equations (LR-Newton-ADI), the Riccati
alternating direction implicit method (RADI)~\cite{BenBKetal18}, the incremental
low-rank subspace iteration (ILRSI)~\cite{LinS15}, and projection-based methods
that construct approximating subspaces such that internally small-scale dense
Riccati equations need to be solved, e.g.,~\cite{HeyJ09, Sim16}.
See also~\cite{BenS13, BenBKetal20, Kue16} for overviews and numerical
comparisons of large-scale sparse solvers for this special case of Riccati
equations with negative semi-definite quadratic terms.
For the opposite case of Riccati equations with positive semi-definite
quadratic terms ($B_{2} = 0$), only the Newton method from~\cite{BenS14}
is known to have a low-rank extension to the large-scale sparse matrix case.
However, there are no methods known to solve the general case~\cref{eqn:indric}
with an indefinite quadratic term in the large-scale sparse setting.
The goal of this paper is to develop an extension of the approach described
in~\cite{LanFA07, LanFAetal08} that can be applied in the case of large-scale
sparse coefficient matrices.

In the upcoming \Cref{sec:rimethod}, we first recap the algorithm
from~\cite{LanFA07, LanFAetal08} and summarize some important results.
Afterwards, we reformulate the steps of the algorithm to fit the large-scale
sparse system case and extend it further to singular $E$ matrices arising in
certain applications.
In \Cref{sec:experiments}, we test the new algorithm first on Riccati equations
with dense coefficients, for a comparison to other established methods, and
afterwards on equations with large-scale sparse coefficients.
The paper is concluded in \Cref{sec:conclusions}.


\section{Riccati iteration method}%
\label{sec:rimethod}

In this section, we describe the idea of the Riccati iteration method
from~\cite{LanFA07, LanFAetal08} and extend the approach to large-scale sparse
systems.

Thereby, we will use the following notation.
We will denote symmetric positive semi-definite matrices by $X = X^{\trans}
\geq 0$, for $X \in \R^{n \times n}$, and use the Loewner ordering $X_{1}
\geq X_{2}$ for two symmetric matrices $X_{1}, X_{2} \in \R^{n \times n}$,
if $X_{1} - X_{2} \geq 0$. 
We call the matrix triple $(A, B, E)$ stabilizable, with $A, E \in
\R^{n \times n}$, $B \in \R^{n \times m}$ and $E$ invertible, if there exists
a feedback matrix $K \in \R^{m \times n}$ such that the matrix pencil
$\lambda E - (A - BK)$ has only eigenvalues with negative real parts.


\subsection{Basic algorithm}

In this section, we formulate the fundamental algorithms for the iterative
computation of approximate solutions to the following general problem.

\begin{problem}[Stabilizing solutions of indefinite Riccati equations]%
  \label{prb:indricsol}
  Given matrices $A, E \in \R^{n \times n}$,  $B_{1} \in \R^{n \times m_{1}}$,
  $B_{2} \in \R^{n \times m_{2}}$ and $C \in \R^{p \times n}$, with $E$
  invertible and  $n$, $m_1$, $m_2$, $p\in \mathbb N$, compute a matrix
  $X = \Xinf\in \R^{n \times n}$, if it exists,
  \begin{itemize}
    \item[(a)] which solves the Riccati equation~\eqref{eqn:indric} with
      indefinite quadratic term
      \begin{align*}
        A^{\trans} X E + E^{\trans} X A + E^{\trans} X \left(
          B_{1} B_{1}^{\trans} - B_{2} B_{2}^{\trans} \right) X E +
          C^{\trans} C & = 0,
        \quad\text{and}
      \end{align*}
    \item[(b)] which is symmetric positive semi-definite and such that the
      pencil
      \begin{align*}
        \lambda E - \big(A + (B_{1} B_{1}^{\trans} - B_{2} B_{2}^{\trans})
          \Xinf E\big)
      \end{align*}
      is stable, i.e., all its eigenvalues lie in the open left half-plane.
  \end{itemize}
\end{problem}

For computing the solution to \Cref{prb:indricsol} numerically, the
\emph{Riccati iteration} method was first described
in~\cite{LanFA07, LanFAetal08} for the standard equation case ($E = I_{n}$).
Here, we will summarize the method and some theoretical results using directly
the generalized case with $E$ invertible.
The underlying idea of the algorithm is to consider the Riccati operator
\begin{align} \label{eqn:ricop}
  \cR(X) & := A^{\trans} X E + E^{\trans} X  A  -
    E^{\trans} X (B_{2} B_{2}^{\trans} - B_{1} B_{1}^{\trans}) X E
    + C_{1}^{\trans} C_{1},
\end{align}
associated with~\cref{eqn:indric}, and a splitting of the final solution
into the sum of consecutive solutions of updated algebraic Riccati equations
with semi-definite quadratic terms.
For two symmetric matrices $X_{1} = X_{1}^{\trans}$ and $X_{2} =
X_{2}^{\trans}$, one can show that
\begin{align*} 
  \cR(X_{1} + X_{2}) & = \cR(X_{1}) + \tA^{\trans} 
    X_{2} E + E^{\trans} X_{2} \tA - E^{\trans} X_{2} 
    (B_{2}B_{2}^{\trans} - B_{1}B_{1}^{\trans}) X_{2} E
\end{align*}
holds, where $\tA := A - (B_{2} B_{2}^{\trans} - 
B_{1} B_{1}^{\trans}) X_{1} E$.
Accordingly, for $X_2$ as a solution to the algebraic
Riccati equation with negative semi-definite quadratic term
\begin{align} \label{eqn:resric}
  0 & = \cR(X_{1}) + \tA^{\trans} 
    X_{2} E + E^{\trans} X_{2} \tA - E^{\trans} X_{2} 
    B_{2}B_{2}^{\trans} X_{2} E,
\end{align}
 the residual reads
\begin{align*}
  \cR(X_{1} + X_{2}) & = E^{\trans} X_{2} B_{1}
    B_{1}^{\trans} X_{2} E.
\end{align*}
An iterative use of this relation~\cref{eqn:resric}, together with
the initial solution $X_{0} = 0$, yields the \emph{Riccati iteration~(RI)}
method.
The resulting algorithm is given in \Cref{alg:ri}.

\begin{algorithm}[t]
  \SetAlgoHangIndent{1pt}
  \DontPrintSemicolon
  \caption{Riccati iteration (RI).}%
  \label{alg:ri}
  
  \SetKw{Error}{error}
  
  \KwIn{$A, B_{1}, B_{2}, C, E$ from~\cref{eqn:indric}, convergence
    tolerance $\tau$.}
  \KwOut{Stabilizing approximate solution $X_{k} \approx \Xinf$
    of~\cref{eqn:indric}, with $X_{k} = X_{k}^{\trans} \geq 0$.}
  
  Initialize $X_{0} = 0$, $k = 0$.\;
  
  \Repeat{$\lVert B_{1} W_{k} E \rVert_{2} < \tau$}{
    Update the iteration matrix $A_{k} = A + (B_{1} B_{1}^{\trans} -
      B_{2} B_{2}^{\trans}) X_{k} E$.\;
      
    Solve the definite Riccati equation
      \begin{align*}
        A_{k}^{\trans} W_{k+1} E + E^{\trans} W_{k+1} A_{k}
          - E^{\trans} W_{k+1} B_{2} B_{2}^{\trans} W_{k+1} E
          + \cR(X_{k}) & = 0
      \end{align*}
    for the stabilizing solution $W_{k+1} \in \R^{n \times n}$.\;
    
    Update the solution matrix $X_{k+1} = X_{k} + W_{k+1}$.\,
    
    \eIf{$(A + B_{1} B_{1}^{\trans} X_{k+1} E, B_{2}, E)$ is stabilizable}{
      Increment $k \leftarrow k + 1$.\;
    }{
      \Error \tcp{There is no stabilizing $X = X^{\trans} \geq 0$
        that solves~\cref{eqn:indric}.}
    }
  }
\end{algorithm}

The following proposition lays out the theoretical foundation of the Riccati
iteration method. 

\begin{proposition}[Properties of the Riccati
  iteration~\cite{LanFA07, LanFAetal08}]%
  \label{prp:ri}
  If $(A, B_{2}, E)$ is stabilizable, $(A, C, E)$ has no unobservable purely 
  imaginary modes and there exists a stabilizing solution
  $\Xinf = \Xinf^{\trans} \geq 0$ for~\cref{eqn:indric}, then the following
  statements hold for the iteration in \Cref{alg:ri}:
  \begin{enumerate}[label=(\alph*)]
    \item $(A + B_{1}B_{1}^{\trans} X_{k}E, B_{2}, E)$ is stabilizable for all
      $k = 0, 1, \ldots$,
    \item $W_{k} = W_{k}^{\trans}  \geq 0$ for all $k = 0, 1, \ldots$,
    \item the eigenvalues of the matrix pencil $\lambda E - (A + 
      B_{1}B_{1}^{\trans} X_{k}E - B_{2} B_{2} X_{k+1} E)$ lie in the left 
      open half-plane for all $k = 0, 1, \ldots$,
    \item $\cR(X_{k+1}) = E^{\trans} W_{k+1} B_{1} B_{1}^{\trans} W_{k+1} E$
      for all $k = 0, 1, \ldots$,
    \item $\Xinf \geq \ldots \geq X_{k+1} \geq \ldots \geq X_{1} \geq
      X_{0} = 0$,
    \item the iteration converges to the stabilizing solution 
      of~\cref{eqn:indric}, $\lim_{k \rightarrow \infty} X_{k} = \Xinf$, and
    \item the convergence is locally quadratic.
  \end{enumerate}
\end{proposition}
Basically, the conditions of \Cref{prp:ri} guarantee convergence of the iterants
to the desired solution of \Cref{prb:indricsol} as well as the stabilization
property for every intermediate iteration step.  
The following remark gives some insight about the convergence of the
stabilizing solutions of~\cref{eqn:indric} and the convergence of the Riccati
iteration.

\begin{remark}[Definiteness of stabilizing solutions]
  In general, if a real stabilizing solution of~\cref{eqn:indric} exists, it
  does not need to be symmetric positive semi-definite; see examples
  in~\cite{LanFA07}.
  Other approaches that rely on the underlying Hamiltonian matrix pencil
  of~\cref{eqn:indric} are capable of computing also indefinite stabilizing
  solutions.
  However, these are undesired in many applications.
  The Riccati iteration converges only if the stabilizing solution exists and
  it is symmetric positive semi-definite, which follows from Parts~(e) and~(f)
  of \Cref{prp:ri} and classical theory about the LQR/LQG Riccati equations
  that are solved in every iteration step.
\end{remark}


\subsection{Factorized low-rank formulation for large-scale sparse equations}

The Riccati iteration method in \Cref{alg:ri} cannot be directly
applied to the large-scale case.
Most importantly, the iterants $X_{k}$ will be dense $n\times n$-matrices such
that memory will become a limiting factor already for moderate dimensions $n$.
With the guarantee that all intermediate solutions $X_{k}$ as well as all
updates $W_{k + 1}$ are symmetric positive semi-definite
(cf. Parts~(a) and~(e) of \Cref{prp:ri}), the approximation by low-rank
factorizations provides a potential remedy as in the symmetric semi-definite
case; see, e.g.,~\cite{BenS13, BenBKetal20}.
The basic idea is to rewrite the intermediate stages of the Riccati iteration
via Cholesky-like low-rank factorizations such that
$X_{k} \approx Z_{k} Z_{k}^{\trans}$ and $W_{k + 1} \approx Y_{k+1}
Y_{k+1}^{\trans}$, with $Z_{k} \in \R^{n \times r_{k}}$ and $Y_{k+1} \in
\R^{n \times q_{k+1}}$ for all $k = 0, 1, \ldots$ and the relevant equations in
terms of the factors $Z_k$ and $Y_k$.
In particular, with the intermediate solution in Step~6 of \Cref{alg:ri} 
reformulated as
\begin{align*}
  Z_{k+1} Z_{k+1}^{\trans} \approx X_{k+1} = X_{k} + W_{k+1} \approx
    Z_{k} Z_{k}^{\trans} + Y_{k+1} Y_{k+1}^{\trans}
    = \begin{bmatrix} Z_{k} & Y_{k+1} \end{bmatrix}
    \begin{bmatrix} Z_{k}^{\trans} \\ Y_{k+1}^{\trans} \end{bmatrix},
\end{align*}
we will develop an iteration for the solution factors without ever
forming the full solution explicitly.

The residual Riccati equation in Step~4 of \Cref{alg:ri} that defines the update $W_{k+1}$ is given via
\begin{align} \label{eqn:resric2}
  A_{k}^{\trans} W_{k+1} E + E^{\trans} W_{k+1} A_{k}
    - E^{\trans} W_{k+1} B_{2} B_{2}^{\trans} W_{k+1} E
    + \cR(X_{k}) & = 0,
\end{align}
and we observe that 
\begin{enumerate}
  \item in the initial step with $X_{0} = 0$, we have
    $\cR(X_{0}) = C^{\trans} C$, and
  \item for all all further steps, with Part~(d) of \Cref{prp:ri}, it holds
    \begin{align*}
      \cR(X_{k+1}) & = E^{\trans} W_{k+1} B_{1} B_{1}^{\trans} W_{k+1} E
        = \underbrace{\big( B_{1}^{\trans} W_{k+1} E \big)^{\trans}}_{
        \phantom{\,\R^{n \times m_{1}}}\in\,\R^{n \times m_{1}}}
        \underbrace{\big( B_{1}^{\trans} W_{k+1} E \big)}_{
        \phantom{\,\R^{m_{1} \times n}}\in\,\R^{m_{1} \times n}},
    \end{align*}
    for all $k = 0, 1, \ldots$, and where
    $W_{k+1} \approx Y_{k+1} Y_{k+1}^{\trans}$.
\end{enumerate}
In both cases, this central step of the algorithm requires the solve of a
standard (semi-definite) Riccati equation with low-rank factorized quadratic
and constant terms.
Accordingly, a low-rank factorized approximation to $X_{k+1}$ can be computed
by established Riccati equation solvers, like the LR-Newton-ADI
method~\cite{BenLP08, Wei16}, RADI~\cite{BenBKetal18},
ILRSI~\cite{LinS15}, or projection-based methods~\cite{HeyJ09, Sim16}.

\begin{remark}[Low-rank solutions]
  By the low-rank structure of the right hand-side and the quadratic terms,
  $m_{1}, m_{2}, p \ll n$, we can also expect the solution and update terms
  to be well approximated by low-rank factors such that $r_{k}, q_{k} \ll n$;
  see~\cite{BenB16, Sti18}.
\end{remark}

Combining all these ideas leads now to the \emph{low-rank Riccati iteration
(LR-RI)} method in \Cref{alg:lrri}.
\Cref{alg:lrri} has been implemented for dense coefficients in the MORLAB
toolbox~\cite{morBenW19b, morBenW21c} and for the large-scale sparse case in the
M-M.E.S.S. toolbox~\cite{SaaKB21-mmess-2.1, morBenKS21}.

\begin{algorithm}[t]
  \SetAlgoHangIndent{1pt}
  \DontPrintSemicolon
  \caption{Low-rank Riccati iteration (LR-RI).}%
  \label{alg:lrri}
  
  \KwIn{$A, B_{1}, B_{2}, C, E$ from~\cref{eqn:indric}, convergence
    tolerance $\tau$.}
  \KwOut{Low-rank factor $Z_{k}$ s.t. $Z_{k}Z_{k}^{\trans} \approx
    \Xinf$ is the stabilizing solution of~\cref{eqn:indric}.}
  
  Solve the Riccati equation with negative semi-definite quadratic term
    \begin{align*}
      A^{\trans} W_{0} E + E^{\trans} W_{0} A 
        - E^{\trans} W_{0} B_{2} B_{2}^{\trans} W_{0} E + C^{\trans} C & = 0,
    \end{align*}
    for the low-rank factor $Y_{0}$ such that
    $Y_{0} Y_{0}^{\trans} \approx W_{0}$.\;
  
  Initialize $Z_{0} = Y_{0}$, $U = \begin{bmatrix} B_{1} & -B_{2}
    \end{bmatrix}$ and $k = 0$.\;
  
  \While{$\lVert B_{1}^{\trans} Y_{k} Y_{k}^{\trans}E \rVert_{2} > \tau$}{
    Set $A_{k} = A - U V_{k}^{\trans}$ with the updated low-rank factor
      \begin{align*}
        V_{k} & = \begin{bmatrix} E^{\trans} Z_{k} Z_{k}^{\trans} B_{1} &
          E^{\trans} Z_{k} Z_{k}^{\trans} B_{2} \end{bmatrix}.
      \end{align*}\vspace{-\baselineskip}\;
    
    Solve the Riccati equation with negative semi-definite quadratic term
      \begingroup \small
      \begin{align*}
        A_{k}^{\trans} W_{k+1} E + E^{\trans} W_{k+1} A_{k}
          - E^{\trans} W_{k+1} B_{2} B_{2}^{\trans} W_{k+1} E
          + \big( E^{\trans} W_{k} B_{1} \big)
          \big( E^{\trans} W_{k} B_{1} \big)^{\trans} & = 0,
      \end{align*}
      \endgroup
      for the low-rank factor $Y_{k+1}$ with
      $Y_{k+1}Y_{k+1}^{\trans} \approx W_{k+1}$.\;
    
    Update $Z_{k+1} = \begin{bmatrix} Z_{k} & Y_{k+1} \end{bmatrix}$ and
      increment $k \leftarrow k + 1$.\;
  }
\end{algorithm}

Another difference of \Cref{alg:lrri} compared to \Cref{alg:ri} is the
stabilizability test for $(A + B_{1} B_{1}^{\trans} X_{k+1} E, B_{2}, E)$.
This additional test is expensive in the large-scale setting but can,
in principle, be omitted.
In case of $(A + B_{1} B_{1}^{\trans} X_{k+1} E, B_{2}, E)$ not being
stabilizable in intermediate steps, the iteration will diverge and
there might be no stabilizing solutions for the intermediate Riccati equations
with negative semi-definite quadratic terms anymore.
The following remarks state further computational features of the new
low-rank iteration method in \Cref{alg:lrri}.

\begin{remark}[Unstable closed-loop matrix pencils]
  The intermediate closed-loop matrix pencils $s E - A_{k}$ can potentially
  have unstable eigenvalues~\cite{LanFAetal08}.
  However, in case there exists a stabilizing solution
  $\Xinf = \Xinf^{\trans} \geq 0$ of~\cref{eqn:indric}, $(A_{k}, B_{2}, E)$
  will be stabilizable.
  Many of the low-rank solvers for definite Riccati equations need a stabilizing
  initial solution or a corresponding feedback matrix.
  In the large-scale sparse case, this can be computed using a sparse eigenvalue
  solver to compute the eigenvectors corresponding to the unstable eigenvalues
  of $s E - A_{k}$ and to apply a partial stabilization approach~\cite{Ben11}
  on the projected problem using the computed eigenvalue
  basis~\cite{BaeBSetal15}.
  Note that if $s E - A_{k_{0}}$ is stable for some $k_{0}$ and the iteration
  converges, then $s E - A_{k}$ will also be stable for all $k \geq k_{0}$;
  cf. \Cref{prp:ri}.
\end{remark}

\begin{remark}[Solution of Riccati equations with positive semi-definite
  quadratic terms]
  The Riccati iteration method can be used as an alternative approach to solve
  Riccati equations with positive semi-definite quadratic terms by setting
  $B_{2} = 0$:
  \begin{align*}
    A^{\trans} X E + E^{\trans} X A + E^{\trans} X B_{1} B_{1}^{\trans} X E
      + C^{\trans} C & = 0,
  \end{align*}
  which occur, for example, in~\cite{morDesP82, morOpdJ88, AndV72}.
  A necessary condition for applying the Riccati iteration method is
  in this case the stability of the matrix pencil $sE - A$, since otherwise the
  stabilizability of $(A, B_{2}, E)$ will never be fulfilled.
  This problem does not occur in the Newton iteration
  from~\cite{BenS14, Var95b, Ben97a}, where only a stabilizing initial feedback
  is needed.
\end{remark}


\subsection{Realization of linear solves in factored form}%
\label{subsec:ak-sparse-lr}

In general, the coefficients $A_{k} = A + (B_{1} B_{1}^{\trans} - B_{2}
B_{2}^{\trans}) X_{k} E$, for $k = 1, 2, \ldots$, of the intermediate Riccati
equations in Step~4 of \Cref{alg:lrri} are $n\times n$-matrices without
sparsity structures such that an explicit realization would make the
approach infeasible in the large-scale setting.
In this case, the coefficients need to be rewritten as the sparse system matrix
$A$ plus low-rank update:
\begin{align}\label{eqn:ak-sparse-plus-lr}
  A_{k} & = A + (B_{1} B_{1}^{\trans} - B_{2} B_{2}^{\trans}) X_{k} E
    = A + \begin{bmatrix} B_{1} & -B_{2} \end{bmatrix}
      \begin{bmatrix} B_{1}^{\trans} X_{k} E \\ B_{2}^{\trans} X_{k} E
      \end{bmatrix}
    =: A + U V_{k}^{\trans},
\end{align}
with $U, V_{k} \in \R^{n \times (m_{1} + m_{2})}$, for all $k = 1, 2, \ldots$.
Thereby, matrix-vector multiplications with $A_{k}$ can be performed without
resorting to an explicit formation of $A_{k}$.

For the use within sparse direct solvers, this factored representation of
$A_{k}$ can be exploited as follows.
Consider the shifted linear system
\begin{align} \label{eqn:linsolve}
  (\sigma E^{\trans} - A_{k}^{\trans}) X &  =
    (\sigma E^{\trans} - A^{\trans} - V_{k} U^{\trans}) X = F,
\end{align}
for a \emph{slim} right-hand side $F \in \R^{n \times \ell}$, with $\ell \ll n$,
and a shift $\sigma \in \C$.
The solution of such linear systems~\cref{eqn:linsolve} are the backbone of
standard iterative solvers for large-scale Riccati equations.
By the \emph{Sherman-Morrison-Woodbury formula} for matrix
inversion~\cite{GolV13} and with the abbreviation
$\Phi(\sigma) := (\sigma E^{\trans} - A^{\trans})$, the inverse of the shifted
matrix in~\cref{eqn:linsolve} is given by
\begin{align} \label{eqn:smw}
  \big( \Phi(\sigma) - V_{k} U^{\trans} \big)^{-1} & =
    \Phi(\sigma)^{-1} + \Phi(\sigma)^{-1} V_{k} \big( I_{m_{1} + m_{2}}
    - U^{\trans} \Phi(\sigma)^{-1} V_{k} \big)^{-1} U^{\trans}
    \Phi(\sigma)^{-1},
\end{align}
which actually amounts to $2 \ell$ sparse linear solves with
$\sigma E^{\trans} - A^{\trans}$ and one solve with an
$(m_{1} + m_{2}) \times (m_1+m_2)$-dimensional matrix.
Thus, in a practical realization of~\cref{eqn:smw} for~\cref{eqn:linsolve}, one
would compute
\begin{align*}
  X & = Z_{1} + Z_{2} \big( I_{m_{1} + m_{2}} - U^{\trans}
    Z_{2} \big)^{-1} U^{\trans} Z_{1},
\end{align*}
with $(\sigma E^{\trans} - A^{\trans}) Z_{1} = F$, and
$(\sigma E^{\trans} - A^{\trans}) Z_{2} = V_{k}$.

An alternative to the Sherman-Morrison-Woodbury formula is the augmented matrix
approach, which makes use of the \emph{block matrix inversion formula}.
This approach is, often times, more stable than the
Sherman-Morrison-Woodbury formula.
Thereby, the solution of~\cref{eqn:linsolve} is given as the solution of
an augmented system of linear equations with
\begin{align} \label{eqn:augmtx}
  \begin{bmatrix} \sigma E^{\trans} - A^{\trans} & V_{k} \\ U^{\trans} &
    I_{m_{1} + m_{2}} \end{bmatrix}
    \begin{bmatrix} X \\ X^{\perp} \end{bmatrix} & =
    \begin{bmatrix} F \\ 0 \end{bmatrix},
\end{align}
where the lower block of the solution, $X^{\perp}$, is an auxiliary variable
with no relevance for the solution.
Under the assumption that $U$ and $V_{k}$ have much fewer columns than $n$,
systems of the form~\cref{eqn:augmtx} can still be solved with standard
sparse solvers.


\subsection{Singular \texorpdfstring{$\boldsymbol{E}$}{E} matrices and
  projected Riccati equations}

A regular occurring case in applications involves singular $E$ matrices.
These especially occur in control problems with differential-algebraic
equations (DAEs).
In general, the presence of singular $E$ matrices changes the concepts of
stabilizability of matrix pencils and solvability of Riccati equations;
see, e.g.,~\cite{morMoeRS11, Hei16, BenH20} and references therein.
However, in many cases, the singular part of $E$ does not play any substantial
role and it is actually enough to consider the solution of Riccati equations
restricted to subspaces corresponding to the finite eigenvalues of the
matrix pencil $\lambda E - A$, i.e., a restriction to the invertible part of
$E$.
This can be realized in two different concepts:
\begin{enumerate}
  \item consider a projected version of the Riccati equation~\cref{eqn:indric}
    with additional constraints on the stabilizing solution,
  \item consider a truncated version of the Riccati equation~\cref{eqn:indric}.
\end{enumerate}

In the following subsections, we consider first the general ideas for projected
and truncated Riccati equations and, afterwards, the setup of incompressible
flows as a particular example.

  
\subsubsection{Projected and truncated equations}%
\label{subsubsec:projric}

In the first concept, appropriate spectral projections $\cP_{r}, \cP_{\ell}
\in \R^{n \times n}$ onto the right and left deflating subspaces corresponding
to the finite eigenvalues of $\lambda E - A$ are explicitly introduced
in~\cref{eqn:indric} such that
\begin{subequations} \label{eqn:projric}
\begin{align}
  A^{\trans} X E + E^{\trans} X A + E^{\trans} X \left( B_{1} B_{1}^{\trans} -
    B_{2} B_{2}^{\trans} \right) X E + \cP_{r}^{\trans} C^{\trans} C \cP_{r}
    & = 0,\\
    \cP_{\ell}^{\trans} X \cP_{\ell} & = X,
\end{align}
\end{subequations}
needs to be solved instead.
While in general, such projections $\cP_{r}, \cP_{\ell}$ can be constructed using
decompositions of the matrix pencil $\lambda E - A$ that resemble the
\emph{Weierstrass canonical form}, they are in fact known for several specially
structured problems that arise in the large-scale sparse setting;
see~\cite{Sty02} and the examples in~\cite{Sty08, morBenS17}.
Note that the $E$ matrix in~\cref{eqn:projric} is still singular, but the
equation and its solutions are restricted to the appropriate underlying subspace.
Iterative methods like \Cref{alg:lrri} can be used to solve~\cref{eqn:projric},
where only the application of the spectral projections $\cP_{r}$ and
$\cP_{\ell}$ is needed in the intermediate computational steps.
This is also the common drawback of solving~\cref{eqn:projric}, since the
repeated application of the projections can quickly become expensive, especially
in the large-scale sparse case.

The second concept of truncated Riccati equations is more commonly used in
practice.
Thereby, we consider in general Riccati equations of the form
\begin{align} \label{eqn:trunric}
  \hA^{\trans} \hX \hE + \hE^{\trans} \hX \hA + \hE^{\trans} \hX \left(
    \hB_{1} \hB_{1}^{\trans} - \hB_{2} \hB_{2}^{\trans} \right) \hX \hE +
    \hC^{\trans} \hC & = 0,
\end{align}
where the truncated matrices are constructed by
\begin{align} \label{eqn:truncoeff}
  \begin{aligned}
    \hA & = W^{\trans} A V, &
      \hE & = W^{\trans} E V, &
      \hB & = W^{\trans} B, &
      \hC & = C V,
  \end{aligned}
\end{align}
with the coefficient matrices from~\cref{eqn:indric} and $V, W \in
\R^{n \times r}$, basis matrices of the right and left deflating subspaces
corresponding to the finite eigenvalues of $\lambda E - A$.
By construction, the $E$ matrix in~\cref{eqn:trunric} is nonsingular and the
solution techniques mentioned so far for~\cref{eqn:indric} can also be
applied to~\cref{eqn:trunric}.
These basis matrices $V$ and $W$ can be computed and applied explicitly
to~\cref{eqn:indric} up to medium-scale sized coefficient matrices.
This is backbone of the implementation of model reduction methods for descriptor
systems in the MORLAB toolbox~\cite{morBenW19b, morBenW21c} and explained, e.g.,
in~\cite{morBenW18}.

However, the explicit computation of the basis matrices $V$ and $W$ and of the
resulting truncated coefficient matrices~\cref{eqn:truncoeff} is usually
not possible in the large-scale sparse setting because of computation time and
memory limitations.
In this case, for many special sparse structures, the truncated Riccati
equation~\cref{eqn:trunric} can be realized implicitly during the computations,
i.e., instead of applying \Cref{alg:lrri} directly to~\cref{eqn:trunric}, the
original sparse coefficients from~\cref{eqn:indric} are used and, by means of
their structure, $V$ and $W$ are implicitly applied during the computational
steps.
Examples for this implicit truncation are given for different sparsity structures
in~\cite{morSaaV18, BaeBSetal15, morHeiSS08, morFreRM08}, which are also
implemented in the function handle framework of the M-M.E.S.S.
toolbox~\cite{SaaKB21-mmess-2.1, morBenKS21}.
Consequently, the implementation of \Cref{alg:lrri} in the M-M.E.S.S. toolbox
can make use of the implicit truncation in case of singular $E$ matrices.
As a particular example, the following subsection considers the case of
structured coefficient matrices arising from incompressible flows.

  
\subsubsection{Implicit realization of truncations in case of flow problems}%
\label{subsubsec:flowsys}

Riccati equations with indefinite quadratic terms~\cref{eqn:indric} are of
particular interest in the design of robust controllers for the stabilization
of incompressible flows modeled by linearization of the Navier-Stokes
equations~\cite{morBenHW21}.
We briefly touch this particular case as an example for implicit truncation
since we will consider such a numerical example later.
In this case, the coefficient matrices of~\cref{eqn:indric} are structured and
given by
\begin{align} \label{eqn:flowsys}
  \begin{aligned}
     A & = \begin{bmatrix} \tA & J^{\trans} \\ J & 0 \end{bmatrix}, &
       E & = \begin{bmatrix} \tE & 0 \\ 0 & 0 \end{bmatrix}, &
       B_{1} & = \begin{bmatrix} \tB_{1} \\ 0 \end{bmatrix}, &
       B_{2} & = \begin{bmatrix} \tB_{2} \\ 0 \end{bmatrix}, &
       C & = \begin{bmatrix} \tC & 0 \end{bmatrix},
  \end{aligned}
\end{align}
with $\tE$ symmetric and invertible.
The key to implicitly truncate and project Riccati equations
with coefficients like~\cref{eqn:flowsys} is the
\emph{discrete Leray projection}
\begin{align} \label{eqn:leray}
  \Pi & = I_{\tn} - \tE^{-1} J^{\trans} (J \tE^{-1} J^{\trans})^{-1} J;
\end{align}
see~\cite{morHeiSS08} for the most general case.
Then, the truncated and projected Riccati equation that needs to be solved has
the form
\begin{equation} \label{eqn:nsetrunric}
  \Pi^{\trans}  \tA^{\trans}\Pi \hX \tE + \tE^{\trans} \hX \Pi^T \tA\Pi  +
    \tE^{\trans} \hX \Pi^{\trans} \left( \tB_{1} \tB_{1}^{\trans} -
    \tB_{2} \tB_{2}^{\trans} \right)\Pi \hX \tE + \Pi^{\trans} \tC^{\trans} \tC
    \Pi  = 0,
\end{equation}
where $\hX = \Pi \hX \Pi^{\trans}$; see~\cite[Sec. 3.1]{morBenHW21}
and~\cite{BenH17b}.

As mentioned in \Cref{subsubsec:projric}, we can neither explicitly compute the
matrices in~\cref{eqn:nsetrunric} nor the projection~\cref{eqn:leray}.
However, instead one can use the original coefficient matrices with their
special structure~\cref{eqn:flowsys}.
As mentioned in \Cref{subsec:ak-sparse-lr}, we practically need to solve
linear systems like~\cref{eqn:linsolve} in every step of the low-rank
Riccati iteration.
If we consider the coefficients of~\cref{eqn:nsetrunric} in this setting,
the low-rank update matrix~\cref{eqn:ak-sparse-plus-lr} has the form
\begin{align*}
  \Pi^{\trans} \tA^{\trans} \Pi + \tE^{\trans} \hX_{k}
    \Pi^{\trans} (\tB_{1} \tB_{1}^{\trans} - \tB_{2} \tB_{2}^{\trans}) \Pi & =
    \Pi^{\trans}(\tA^{\trans} + \tE^{\trans} \hX_{k} (\tB_{1} \tB_{1}^{\trans} -
    \tB_{2} \tB_{2}^{\trans})) \Pi
    = \Pi^{\trans} \tA_{k} \Pi,
\end{align*}
where we also used that $\tE \Pi = \Pi^{\trans} \tE$ and $\Pi^{2} =\Pi$.
By observing that all occurring right-hand sides will satisfy $F = \Pi F$,
the solution of $(\sigma \tE^{\trans} - \Pi^{\trans} \tA_{k}^{\trans} \Pi) X
= F$ is alternatively given by
\begin{align} \label{eqn:linsolv2}
  \begin{bmatrix} \sigma \tE^{\trans} - \tA^{\trans} + \tE^{\trans} \hX_{k}
    (\tB_{1} \tB_{1}^{\trans} - \tB_{2} \tB_{2}^{\trans}) & -J^{\trans} \\
    -J & 0 \end{bmatrix} \begin{bmatrix} X \\ X^{\perp} \end{bmatrix}
    & = \begin{bmatrix} F \\ 0 \end{bmatrix},
\end{align}
where $X^{\perp}$ is an auxiliary variable that does not play any role.
This ensures that $X = \Pi X$ as required during the iteration;
see~\cite[Lem. 5.2]{morHeiSS08}; and respects the sparsity structure
of~\cref{eqn:flowsys}, since~\cref{eqn:linsolv2} is actually given by
\begin{align*}
  \left( \sigma E^{\trans} - A^{\trans} + E^{\trans}
    \begin{bmatrix} \hX_{k} & 0 \\ 0 & 0 \end{bmatrix}
    (B_{1} B_{1}^{\trans} - B_{2} B_{2}^{\trans}) \right)
    \begin{bmatrix} X \\ X^{\perp} \end{bmatrix}
    & = \begin{bmatrix} F \\ 0 \end{bmatrix},
\end{align*}
using the matrices from~\cref{eqn:flowsys}.


\section{Numerical experiments}%
\label{sec:experiments}

The numerical experiments reported in this section have been executed on a
machine with 2 Intel(R) Xeon(R) Silver 4110 CPU processors running at 2.10\,GHz
and equipped with 192\,GB total main memory.
The computer runs on CentOS Linux release 7.5.1804 (Core) with
MATLAB 9.9.0.1467703 (R2020b).

\begin{table}[t]
  \caption{Results for the aircraft and cable mass benchmarks.}
  \label{tab:examples}
   
  \centering \small
  \renewcommand{\arraystretch}{1.25}
  \begin{tabular}{|rrrrrrrr|}
    \hline \Tstrut \Bstrut
    & 
      $n$ &
      $m_{1}$ &
      $m_{2}$ &
      $p$ &
      \#$\lambda_{\mathrm{unstab}}$ &
      data &
      $\gamma$ \\ \hline \Tstrut
    \texttt{AC10} & $55$ & $3$ & $2$ & $5$ & $2$ & dense & $2.00$ \\
    \texttt{CM6} & $960$ & $1$ & $1$ & $3$ & $76$ & dense & $2.00$ \\
    \texttt{rand512} & $512$ & $2$ & $6$ & $5$ & $8$ & dense & $8.00$ \\
    \texttt{rand1024} & $1\,024$ & $4$ & $2$ & $3$ & $6$ & dense & $8.00$ \\
    \texttt{rand2048} & $2\,048$ & $3$ & $3$ & $1$ & $2$ & dense & $7.00$ \\ 
    \Bstrut
    \texttt{rand4096} & $4\,096$ & $4$ & $3$ & $6$ & $1$ & dense & $10.00$\\
    \hline \Tstrut
    \texttt{rail} & $79\,841$ & $3$ & $3$ & $6$ & $0$ & sparse & $2.00$ \\ 
    \Bstrut
    \texttt{cylinderwake} & $47\,136$ & $1$ & $1$ & $6$ & $4$ & sparse &
      $50.00$ \\
    \hline
  \end{tabular}
\end{table}

The benchmark data used can be found in \Cref{tab:examples}, where the first
four columns state the dimensions of the considered problem;
cf.~\cref{eqn:indric}; \#$\lambda_{\mathrm{unstab}}$ denotes the number of
unstable eigenvalues of the matrix pencil $s E - A$, in the data column the
example type is denoted with either dense or sparse, and $\gamma$ is a constant
multiplied with the disturbance term, $\frac{1}{\gamma} B_{1}$, such that the
considered Riccati equations have symmetric positive semi-definite stabilizing
solutions.

For the comparison of results between different algorithms and benchmarks, we
will report different information about the performance of the applied solvers
in upcoming tables such as the used iteration steps, the overall runtime in
seconds, the rank of the resulting solution factor, the final normalized
residual of the Riccati iteration as given by $\lVert B_{1} Y_{k}
Y_{k}^{\trans} E \rVert_{2}^{2} / \lVert C C^{\trans} \rVert_{2}$,
the relative residual of the solution factor given by $\lVert \cR(Z_{k}
Z_{k}^{\trans}) \rVert_{2} / \lVert Z_{k}^{\trans} Z_{k} \rVert_{2}$ and
the normalized residual of the solution factor computed as $\lVert \cR(Z_{k}
Z_{k}^{\trans}) \rVert_{2} / \lVert C C^{\trans} \rVert_{2}$ using the Riccati
operator $\cR$ from~\cref{eqn:ricop}, and lastly the norm of the full solution
$\lVert Z_{k}^{\trans} Z_{k} \rVert_{2}$.

\begin{center} \small
  \setlength{\fboxsep}{5pt}%
  \fbox{%
    \begin{minipage}{.92\textwidth}
      \textbf{Code and data availability}\newline
      The source code of the implementations used to compute the presented
      results, the used data and computed results are available at
      \begin{center} 
        \href{https://doi.org/10.5281/zenodo.5155993}%
          {\texttt{doi:10.5281/zenodo.5155993}}
      \end{center}
      under the BSD-2-Clause license and authored by Steffen W. R. Werner.
    \end{minipage}%
  }%
\end{center}


\subsection{Factorized method for dense examples}

Before we actually consider the large-scale sparse case, we test the new
method from \Cref{alg:lrri} on medium-scale dense benchmarks and compare the
results with two other commonly known solvers that can be applied
to~\cref{eqn:indric}.
We use the dense low-rank factor version of \Cref{alg:lrri} as implemented
in \texttt{ml\_icare\_ric\_fac} of~\cite{morBenW19b}, which we denote further
on as \textsf{LRRI}, but we replace the internal Riccati equation solver by the
factorized sign function approach~\cite{BenEQetal14} implemented in the very
same toolbox~\cite{morBenW19b} as \texttt{ml\_caredl\_sgn\_fac}.
For comparison, we use the Hamiltonian eigenvalue approach from~\cite{ArnL84}
as it is implemented in the new MATLAB function \texttt{icare} from the
Control System Toolbox\texttrademark{}, further denoted by \textsf{ICARE}, and
the classical (unfactorized) sign function solver~\cite{morRob80} from the
MORLAB toolbox~\cite{morBenW19b} in \texttt{ml\_caredl\_sgn}, further on as
\textsf{SIGN}.
Since \textsf{SIGN} and \textsf{ICARE} can only compute unfactorized solutions
of~\cref{eqn:indric}, we perform eigendecompisitions of the computed solutions
to check the positive semi-definiteness and to compute low-rank approximations
of the solutions by truncating all components corresponding to non-positive
eigenvalues.

Considering the benchmarks, the first two data sets in \Cref{tab:examples}
are practical examples taken from~\cite{Lei04}, where \texttt{AC10} is an
aircraft model and \texttt{CM6} a cable mass model with low damping.
The data was taken over with exactly the same naming as given in~\cite{Lei04},
where we have set $C = C_{1}$.
The results of the different algorithms on these two data sets are shown in
\Cref{tab:accm}.
First, we recognize that for \texttt{AC10}, \textsf{LRRI} performs visibly worse
than the other two approaches in terms of accuracy, since we are loosing three
orders of magnitude in the relative and normalized residuals.
We assume this comes from the bad conditioning of the $A$ matrix in this example
and the repeated solution of Riccati equations with the matrix.
However, for the \texttt{CM6} example, this turned around, as \textsf{LRRI}
performs now an order of magnitude better than \textsf{ICARE} and two orders
better than \textsf{SIGN}.
Concerning the computation times we see that \textsf{LRRI} performs reasonably
well thanks to the efficient inner Riccati equation solver.
While on the smaller data set \texttt{AC10}, the runtime of \textsf{LRRI} is in
the same order of magnitude as \textsf{ICARE}, \textsf{LRRI} is already four
times faster than \textsf{ICARE} for \texttt{CM6}.
We are not able to outperform \textsf{SIGN} with \textsf{LRRI} in these
examples.

\begin{table}[t]
  \caption{Results for the aircraft and cable mass benchmarks.}
  \label{tab:accm}
   
  \centering \small
  \renewcommand{\arraystretch}{1.25}
  \begin{tabular}{|r|rrr|rrr|}
    \hline \Tstrut 
    & \multicolumn{3}{c|}{\texttt{AC10}} &
      \multicolumn{3}{c|}{\texttt{CM6}} \\ \Bstrut
    & \multicolumn{1}{c}{\textsf{LRRI}} &
      \multicolumn{1}{c}{\textsf{ICARE}} &
      \multicolumn{1}{c|}{\textsf{SIGN}} &
      \multicolumn{1}{c}{\textsf{LRRI}} &
      \multicolumn{1}{c}{\textsf{ICARE}} &
      \multicolumn{1}{c|}{\textsf{SIGN}} \\ \hline \Tstrut
    Iteration steps &
      $5$ &
      --- &
      $19$ &
      $4$ &
      --- &
      $23$ \\
    Runtime (s) &
      $0.92801$ &
      $0.47971$ &
      $0.07705$ &
      $29.9208$ &
      $140.136$ &
      $7.24722$ \\
    Rank $Z_{k}$ & 
      $53$ &
      $55$ &
      $55$ &
      $569$ &
      $758$ &
      $781$ \\
    Final res. &
      \mnum{5.545}{-}{25} &
      --- &
      --- &
      \mnum{1.873}{-}{15} &
      --- &
      --- \\
    Relative res. &
      \mnum{2.599}{-}{07} &
      \mnum{9.617}{-}{10} &
      \mnum{1.183}{-}{09} &
      \mnum{1.910}{-}{09} &
      \mnum{5.019}{-}{08} &
      \mnum{2.125}{-}{07}\\
    Normalized res. &
      \mnum{1.554}{-}{03} &
      \mnum{5.752}{-}{06} &
      \mnum{7.074}{-}{06} &
      \mnum{1.667}{-}{05} &
      \mnum{4.381}{-}{04} &
      \mnum{1.855}{-}{03} \\ \Bstrut
    $\lVert Z_{k}^{\trans} Z_{k} \rVert_{2}$ &
      \mnum{1.457}{+}{01} &
      \mnum{1.457}{+}{01} &
      \mnum{1.457}{+}{01} &
      \mnum{1.253}{+}{04} &
      \mnum{1.253}{+}{04} &
      \mnum{1.253}{+}{04} \\[.6ex]
    \hline
  \end{tabular}
\end{table}

\begin{table}[t]
  \caption{Results for the smaller random examples.}
  \label{tab:randsmall}
   
  \centering \small
  \renewcommand{\arraystretch}{1.25}
  \begin{tabular}{|r|rrr|rrr|}
    \hline \Tstrut 
    & \multicolumn{3}{c|}{\texttt{rand512}} &
      \multicolumn{3}{c|}{\texttt{rand1024}} \\ \Bstrut
    & \multicolumn{1}{c}{\textsf{LRRI}} &
      \multicolumn{1}{c}{\textsf{ICARE}} &
      \multicolumn{1}{c|}{\textsf{SIGN}} &
      \multicolumn{1}{c}{\textsf{LRRI}} &
      \multicolumn{1}{c}{\textsf{ICARE}} &
      \multicolumn{1}{c|}{\textsf{SIGN}} \\ \hline \Tstrut
    Iteration steps &
      $8$ &
      --- &
      $14$ &
      $7$ &
      --- &
      $8$ \\
    Runtime (s) &
      $6.94609$ &
      $23.0720$ &
      $1.39393$ &
      $10.3911$ &
      $217.230$ &
      $2.96977$ \\
    Rank $Z_{k}$ & 
      $96$ &
      $298$ &
      $294$ &
      $89$ &
      $522$ &
      $512$ \\
    Final res. &
      \mnum{7.889}{-}{23} &
      --- &
      --- &
      \mnum{2.600}{-}{14} &
      --- &
      --- \\
    Relative res. &
      \mnum{7.544}{-}{10} &
      \mnum{2.926}{-}{11} &
      \mnum{3.508}{-}{11} &
      \mnum{6.737}{-}{11} &
      \mnum{6.980}{-}{11} &
      \mnum{1.074}{-}{10}\\
    Normalized res. &
      \mnum{3.930}{-}{10} &
      \mnum{1.524}{-}{11} &
      \mnum{1.827}{-}{11} &
      \mnum{6.164}{-}{12} &
      \mnum{6.387}{-}{12} &
      \mnum{9.830}{-}{12} \\ \Bstrut
    $\lVert Z_{k}^{\trans} Z_{k} \rVert_{2}$ &
      \mnum{3.202}{+}{02} &
      \mnum{3.202}{+}{02} &
      \mnum{3.202}{+}{02} &
      \mnum{1.058}{+}{02} &
      \mnum{1.058}{+}{02} &
      \mnum{1.058}{+}{02} \\[.6ex]
    \hline
  \end{tabular}
\end{table}

To further investigate accuracy and performance of \textsf{LRRI}, we created
random examples denoted by \texttt{rand*} in \Cref{tab:examples} using the
\texttt{randn} function in MATLAB.
The results of the computations can be found in
\Cref{tab:randsmall,tab:randlarge}.
Overall, \textsf{LRRI} competes very well against \textsf{ICARE} and
\textsf{SIGN} in terms of accuracy.
For all four presented random examples, the residuals of the solution factors
lie in the same order of magnitude for all methods with only minor exceptions.
In terms of runtimes, we see the same relations that we already recognized in
\Cref{tab:accm}.
\textsf{LRRI} easily outperforms \textsf{ICARE} due to its efficient inner
factorized sign function solver.
For increasing problem size also this speed-up tremendously increases further up
to the largest example \texttt{rand4096}, where \textsf{LRRI} is $113$ times
faster than \textsf{ICARE}.
Compared to \textsf{SIGN}, \textsf{LRRI} is still not able to outperform the
classical sign function iteration method.
However, we can observe that the difference in runtimes is getting smaller and 
smaller, i.e., we can expect \textsf{LRRI} to outperform \textsf{SIGN} for
larger problems.
Also, these results reveal the strong dependence of \textsf{LRRI} on the inner
Riccati equation solver and the power that comes from the ability to use
sophisticated implementations for the inner Riccati equations with negative
semi-definite quadratic terms.

\begin{table}[t]
  \caption{Results for the larger random examples.}
  \label{tab:randlarge}
   
  \centering \small
  \renewcommand{\arraystretch}{1.25}
  \begin{tabular}{|r|rrr|rrr|}
    \hline \Tstrut 
    & \multicolumn{3}{c|}{\texttt{rand2048}} &
      \multicolumn{3}{c|}{\texttt{rand4096}} \\ \Bstrut
    & \multicolumn{1}{c}{\textsf{LRRI}} &
      \multicolumn{1}{c}{\textsf{ICARE}} &
      \multicolumn{1}{c|}{\textsf{SIGN}} &
      \multicolumn{1}{c}{\textsf{LRRI}} &
      \multicolumn{1}{c}{\textsf{ICARE}} &
      \multicolumn{1}{c|}{\textsf{SIGN}} \\ \hline \Tstrut
    Iteration steps &
      $5$ &
      --- &
      $11$ &
      $4$ &
      --- &
      $9$ \\
    Runtime (s) &
      $34.8321$ &
      $2043.78$ &
      $21.4784$ &
      $139.670$ &
      $16042.8$ &
      $125.835$ \\
    Rank $Z_{k}$ & 
      $59$ &
      $992$ &
      $967$ &
      $72$ &
      $1984$ &
      $1961$ \\
    Final res. &
      \mnum{2.450}{-}{21} &
      --- &
      --- &
      \mnum{2.000}{-}{14} &
      --- &
      --- \\
    Relative res. &
      \mnum{4.997}{-}{10} &
      \mnum{2.977}{-}{10} &
      \mnum{1.628}{-}{10} &
      \mnum{1.811}{-}{10} &
      \mnum{1.804}{-}{09} &
      \mnum{1.813}{-}{10}\\
    Normalized res. &
      \mnum{3.402}{-}{11} &
      \mnum{2.027}{-}{11} &
      \mnum{1.108}{-}{11} &
      \mnum{2.812}{-}{11} &
      \mnum{2.803}{-}{10} &
      \mnum{2.816}{-}{11} \\ \Bstrut
    $\lVert Z_{k}^{\trans} Z_{k} \rVert_{2}$ &
      \mnum{1.642}{+}{02} &
      \mnum{1.642}{+}{02} &
      \mnum{1.642}{+}{02} &
      \mnum{6.760}{+}{02} &
      \mnum{6.760}{+}{02} &
      \mnum{6.760}{+}{02} \\[.6ex]
    \hline
  \end{tabular}
\end{table}


\subsection{Low-rank approach for large-scale sparse examples}

Now we come to the case of large-scale sparse Riccati equations with indefinite
quadratic terms.
The \textsf{LRRI} method from \Cref{alg:lrri} is to our knowledge the only
method suited to solve equations like~\cref{eqn:indric} in the large-scale
sparse setting.
Therefor, we cannot compare our results to other methods.
However, from the previous section we expect \textsf{LRRI} to yield reasonable
accurate results in comparison to alternative Riccati equation solvers and the
runtimes to be mainly depend on the inner Riccati equation solver.
We use the implementation of \textsf{LRRI} from~\cite{SaaKB21-mmess-2.1} with
the option to have RADI~\cite{BenBKetal18} or LR-Newton-ADI~\cite{BenLP08} as
inner solvers or to switch between them during runtime if suitable.

\begin{figure}[t]
  \centering
  \tikzexternalenable%
  \tikzsetnextfilename{rail}%
  \filemodCmp{graphics/rail.tikz}{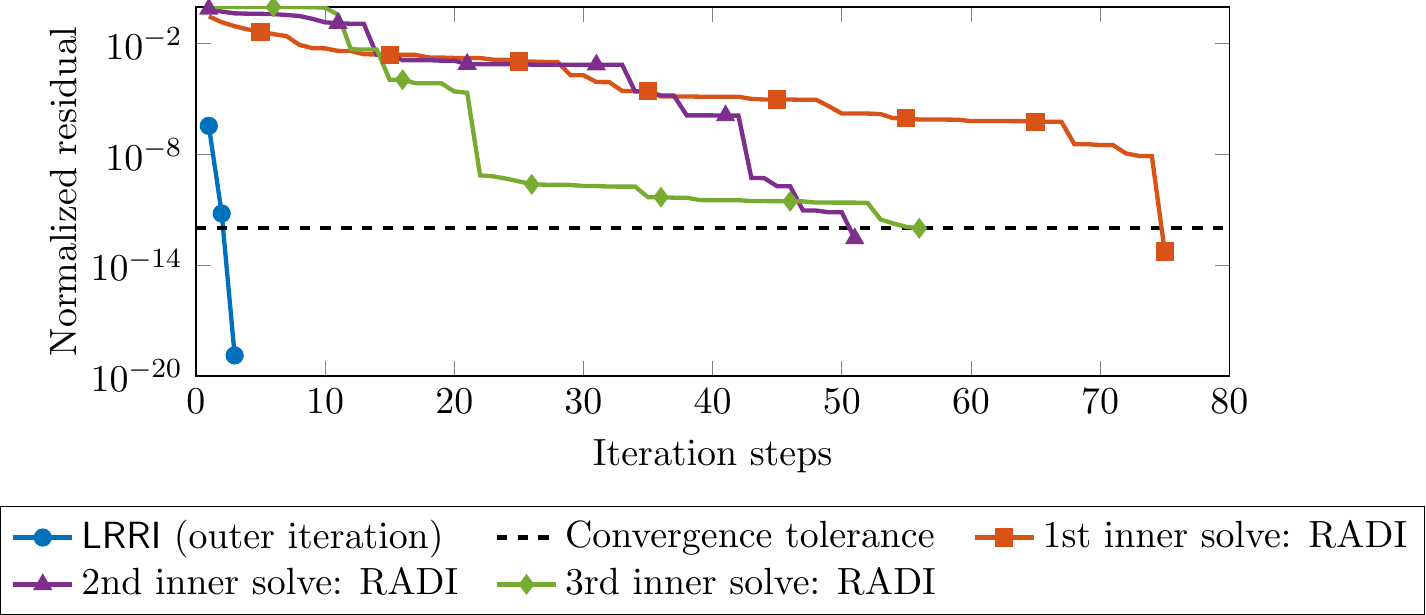}%
    {\tikzset{external/remake next}}{}%
  \begin{tikzpicture}
  \pgfplotstableread{graphics/data/rail_ri.dat}\tableRI
  \pgfplotstableread{graphics/data/rail_radi1.dat}\tableRADIA
  \pgfplotstableread{graphics/data/rail_radi2.dat}\tableRADIB
  \pgfplotstableread{graphics/data/rail_radi3.dat}\tableRADIC
  
  \begin{axis}[%
    width  = .7\textwidth,
    height = .25\textwidth,
    scale only axis,
    xmin = 0,
    xmax = 80,
    ymin = 1e-20,
    ymax = 1e+00,
    ymode = log,
    xminorticks = false,
    yminorticks = false,
    xlabel = {Iteration steps},
    ylabel = {Normalized residual},
    ylabel style   = {yshift = -.3em},
    scaled x ticks = false,
    x tick label style = {/pgf/number format/fixed},
    legend columns = 3,
    legend style = {
      at     = {(.5,-.5)},
      anchor = center,
      /tikz/every even column/.append style = {column sep = 0.3cm}},
    legend cell align = {left}]
      
    \addplot[ri] table[x index = 0, y index = 1] {\tableRI};
    \addlegendentry{\textsf{LRRI} (outer iteration)}
      
    \addplot[tol] coordinates {(0,1e-12) (80, 1e-12)};
    \addlegendentry{Convergence tolerance}
    
    \addplot[inner1, mark phase = 5, mark repeat = 10]
      table[x index = 0, y index = 1] {\tableRADIA};
    \addlegendentry{1st inner solve: RADI}
    
    \addplot[inner2, mark phase = 1, mark repeat = 10]
      table[x index = 0, y index = 1] {\tableRADIB};
    \addlegendentry{2nd inner solve: RADI}
    
    \addplot[inner3, mark phase = 6, mark repeat = 10]
      table[x index = 0, y index = 1] {\tableRADIC};
    \addlegendentry{3rd inner solve: RADI}
  \end{axis}
\end{tikzpicture}%
  \tikzexternaldisable%

  \caption{Convergence of Riccati iteration and inner solvers for the
    \texttt{rail} example.}
  \label{fig:rail}
\end{figure}
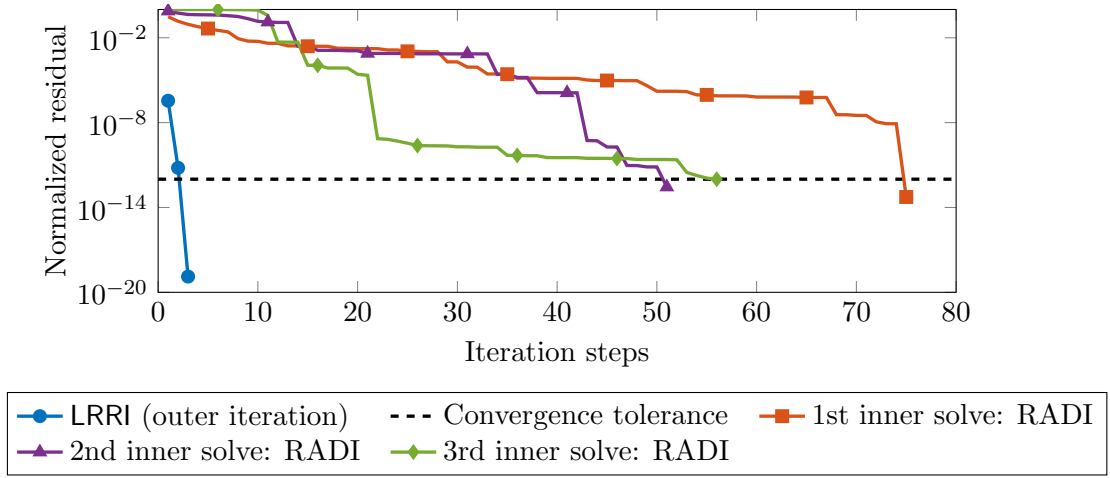

As first example, we consider the optimal cooling problem of a steel
profile; see, e.g.,~\cite{Saa09}; with the data available
in~\cite{SaaKB21-mmess-2.1}.
For the different matrices in the quadratic term of~\cref{eqn:indric}, we
consider the boundary control of the three lower segments of the profile
edges as disturbances to give us $B_{1}$ and the rest to be control inputs
in $B_{2}$.
The resulting dimensions are given in \Cref{tab:examples}.
For \textsf{LRRI}, we use only the RADI method as inner Riccati equation
solver and the results of the iteration can be seen in the \texttt{rail} column
of \Cref{tab:sparse}.
The iteration converges quickly and yields a very accurate solution of
low rank.
The behavior of the outer and inner iteration methods is shown in terms of
normalized residuals in \Cref{fig:rail}.
For all solvers, the same convergence tolerance has been used.

\begin{table}[t]
  \caption{Results of \textsf{LRRI} for the large-scale sparse examples.}
  \label{tab:sparse}
   
  \centering \small
  \renewcommand{\arraystretch}{1.25}
  \begin{tabular}{|r|r|r|}
    \hline \Tstrut \Bstrut
    & \multicolumn{1}{c|}{\texttt{rail}} &
      \multicolumn{1}{c|}{\texttt{cylinderwake}} \\ \hline \Tstrut
    Iteration steps &
      $3$ &
      $3$ \\
    Runtime (s) &
      $72.7378$ &
      $3469.59$ \\
    Rank $Z_{k}$ &
      $169$ &
      $418$ \\
    Final res. &
      \mnum{1.297}{-}{19} &
      \mnum{2.184}{-}{21} \\
    Relative res. &
      \mnum{2.125}{-}{21} &
      \mnum{1.996}{-}{14} \\
    Normalized res. &
      \mnum{9.766}{-}{11} &
      \mnum{1.622}{-}{03} \\ \Bstrut
    $\lVert Z_{k}^{\trans} Z_{k} \rVert_{2}$ &
      \mnum{6.866}{+}{11} &
      \mnum{5.056}{+}{08} \\[.6ex]
    \hline
  \end{tabular}
\end{table}

\begin{figure}[t]
  \centering
  \tikzexternalenable%
  \tikzsetnextfilename{cylinderwake}%
  \filemodCmp{graphics/cylinderwake.tikz}{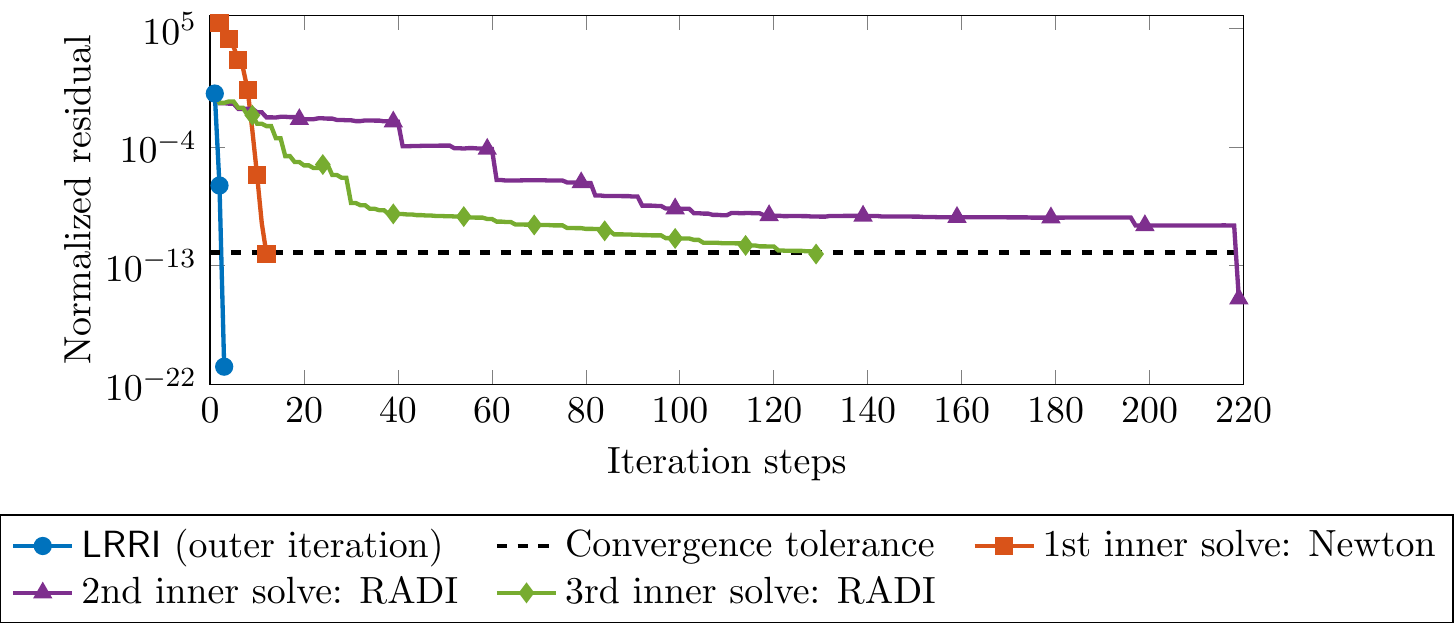}%
    {\tikzset{external/remake next}}{}%
  \begin{tikzpicture}
  \pgfplotstableread{graphics/data/cylinderwake_ri.dat}\tableRI
  \pgfplotstableread{graphics/data/cylinderwake_nwt1.dat}\tableNWTA
  \pgfplotstableread{graphics/data/cylinderwake_radi2.dat}\tableRADIB
  \pgfplotstableread{graphics/data/cylinderwake_radi3.dat}\tableRADIC
  
  \begin{axis}[%
    width  = .7\textwidth,
    height = .25\textwidth,
    scale only axis,
    xmin = 0,
    xmax = 220,
    ymin = 1e-22,
    ymax = 1e+06,
    ymode = log,
    xminorticks = false,
    yminorticks = false,
    xlabel = {Iteration steps},
    ylabel = {Normalized residual},
    ylabel style   = {yshift = -.3em},
    scaled x ticks = false,
    x tick label style = {/pgf/number format/fixed},
    legend columns = 3,
    legend style = {
      at     = {(.5,-.5)},
      anchor = center,
      /tikz/every even column/.append style = {column sep = 0.3cm}},
    legend cell align = {left}]
      
    \addplot[ri] table[x index = 0, y index = 1] {\tableRI};
    \addlegendentry{\textsf{LRRI} (outer iteration)}
      
    \addplot[tol] coordinates {(0,1e-12) (220, 1e-12)};
    \addlegendentry{Convergence tolerance}
    
    \addplot[inner1, mark phase = 0, mark repeat = 2]
      table[x index = 0, y index = 1] {\tableNWTA};
    \addlegendentry{1st inner solve: Newton}
    
    \addplot[inner2, mark phase = 19, mark repeat = 20]
      table[x index = 0, y index = 1] {\tableRADIB};
    \addlegendentry{2nd inner solve: RADI}
    
    \addplot[inner3, mark phase = 9, mark repeat = 15]
      table[x index = 0, y index = 1] {\tableRADIC};
    \addlegendentry{3rd inner solve: RADI}
  \end{axis}
\end{tikzpicture}%
  \tikzexternaldisable%

  \caption{Convergence of Riccati iteration and inner solvers for the
    cylinder wake example.}
  \label{fig:cylinderwake}
\end{figure}
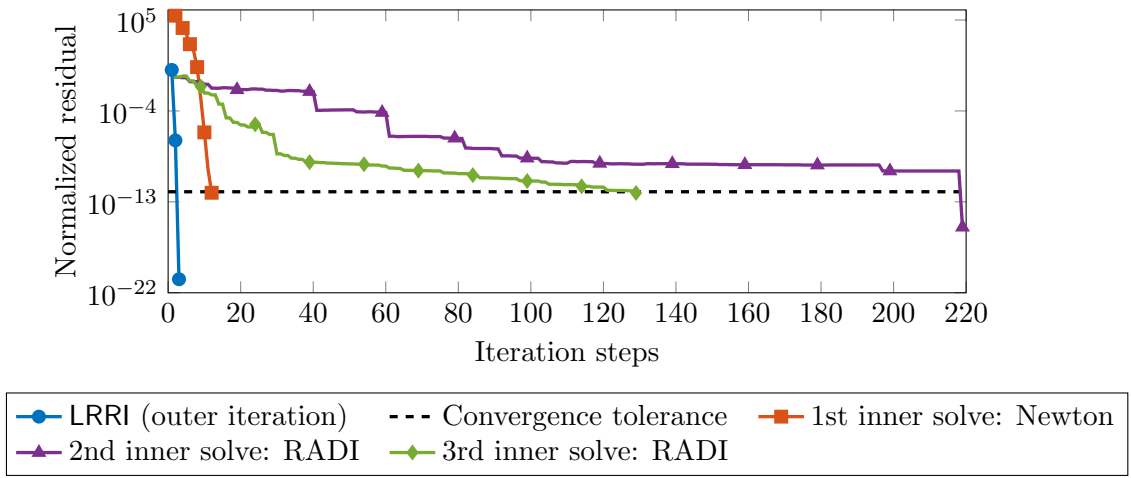

As second example, we consider the laminar flow in a wake with a cylinder
obstacle at Reynolds number $60$, as described in~\cite{morBenHW21} and with
the data available at~\cite{supHeiW21a}.
The example has two inputs modeling the suction and injection of fluid at the
back of the cylinder placed in the beginning of the wake.
We consider the case that one of the outlets is defect and produces only noise
that needs to be compensated.
Therefor, the defect outlet gives us the $B_{1}$ matrix and the control outlet
gives us the matrix $B_{2}$.
Also, the considered Riccati equation has exactly the structure as
in~\cref{eqn:flowsys} resulting from the underlying dynamical system of
differential-algebraic equations.
Therefor, we use the implicit truncation approach mentioned in
\Cref{subsubsec:flowsys}, which is implemented in the \texttt{dae\_2} function
handles in~\cite{SaaKB21-mmess-2.1}.
Also, we note that the matrix pencil of this example has unstable eigenvalues.
Since RADI is difficult to use for such a problem since it needs a stabilizing
initial solution that produces a positive semi-definite residual in the Riccati
operator, we switch only for the first iteration step of \textsf{LRRI} to the
LR-Newton-ADI method and use a stabilizing Bernoulli feedback the same way as
in~\cite{BaeBSetal15}.
The results of \textsf{LRRI} can be seen in the \texttt{cylinderwake} column of
\Cref{tab:sparse} and the normalized residuals of the outer and inner iterations
are shown in \Cref{fig:cylinderwake}.
While the relative residual is comfortably small again for this example due 
to the very large norm of the stabilizing solution, the normalized residual
is still quite high.
This likely comes from the general bad conditioning of the problem and the
small norm of the right-hand side matrix.
But overall, the computed solution has been obtained with reasonable accuracy
and in a reasonable amount of time.

\section{Conclusions}%
\label{sec:conclusions}

We have developed a low-rank iterative method for solving large-scale sparse
Riccati equations with indefinite quadratic terms, which is based on solutions
of Riccati equations with negative semi-definite quadratic terms.
Numerical examples have illustrated that, in the dense case, we can expect a
similar accuracy and good performance in comparison to other established Riccati
equation solvers, and that for large-scale sparse equations, the method also
yields reasonably good results.
We have also extended the LR-RI approach to indefinite algebraic Riccati
equations related to descriptor systems with singular $E$ matrix.

To our knowledge, the low-rank Riccati iteration is currently the only approach
to solve Riccati equations with indefinite quadratic terms in the large-scale
sparse case.
Another idea that might be directly extendable to this problem are
projection-based methods.
However, a problem occurring already in the case of classical Riccati equations
with semi-definite quadratic terms is the solvability of the projected equations.
Only recently, new results have been obtained under which conditions
constructed projection spaces preserve the existence of stabilizing solutions
in projected Riccati equations with negative semi-definite quadratic
terms~\cite{ZhaFC20}.
This problem becomes even more complicated when dealing with Riccati equations
with indefinite quadratic terms, and so far, there are neither theoretical nor 
numerical results available considering this.


\addcontentsline{toc}{section}{Acknowledgment}
\section*{Acknowledgment}

All authors have been supported by the German Research Foundation (DFG) Research
Training Group 2297 ``MathCoRe'', Magdeburg.


\addcontentsline{toc}{section}{References}

\bibliographystyle{plainurl}
\bibliography{bibtex/myref}

\begin{thebibliography}{10}

\bibitem{AmmBM93}
G.~S. Ammar, P.~Benner, and V.~Mehrmann.
\newblock A multishift algorithm for the numerical solution of algebraic
  {R}iccati equations.
\newblock {\em Electron. Trans. Numer. Anal.}, 1:33--48, 1993.
\newblock URL:
  \url{https://etna.math.kent.edu/volumes/1993-2000/vol1/abstract.php?vol=1&pages=33-48}.

\bibitem{AndM90}
B.~D.~O. Anderson and J.~B. Moore.
\newblock {\em Optimal Control: {L}inear Quadratic Methods}.
\newblock Prentice-Hall, Englewood Cliffs, NJ, 1990.

\bibitem{AndV72}
B.~D.~O. Anderson and B.~Vongpanitlerd.
\newblock {\em Network Analysis and Synthesis: {A} Modern Systems Approach}.
\newblock Networks Series. Prentice-Hall, Englewood Cliffs, NJ, 1972.

\bibitem{ArnL84}
W.~F. Arnold and A.~J. Laub.
\newblock Generalized eigenproblem algorithms and software for algebraic
  {R}iccati equations.
\newblock {\em Proc. {IEEE}}, 72(12):1746--1754, 1984.
\newblock \href {https://doi.org/10.1109/PROC.1984.13083}
  {\path{doi:10.1109/PROC.1984.13083}}.

\bibitem{BaeBSetal15}
E.~B{\"a}nsch, P.~Benner, J.~Saak, and H.~K. Weichelt.
\newblock {R}iccati-based boundary feedback stabilization of incompressible
  {N}avier-{S}tokes flows.
\newblock {\em {SIAM} J. Sci. Comput.}, 37(2):A832--A858, 2015.
\newblock \href {https://doi.org/10.1137/140980016}
  {\path{doi:10.1137/140980016}}.

\bibitem{BasM17}
T.~Ba{\c{s}}ar and J.~Moon.
\newblock {R}iccati equations in {N}ash and {S}tackelberg differential and
  dynamic games.
\newblock {\em IFAC-PapersOnLine}, 50(1):9547--9554, 2017.
\newblock \href {https://doi.org/10.1016/j.ifacol.2017.08.1625}
  {\path{doi:10.1016/j.ifacol.2017.08.1625}}.

\bibitem{Ben97}
P.~Benner.
\newblock {\em Contributions to the Numerical Solution of Algebraic {R}iccati
  Equations and Related Eigenvalue Problems}.
\newblock {Dissertation}, Department of Mathematics, TU Chemnitz-Zwickau,
  Chemnitz, Germany, 1997.

\bibitem{Ben97a}
P.~Benner.
\newblock Numerical solution of special algebraic {R}iccati equations via an
  exact line search method.
\newblock In {\em 1997 European Control Conference (ECC)}, pages 3136--3141,
  1997.
\newblock \href {https://doi.org/10.23919/ECC.1997.7082591}
  {\path{doi:10.23919/ECC.1997.7082591}}.

\bibitem{Ben11}
P.~Benner.
\newblock Partial stabilization of descriptor systems using spectral
  projectors.
\newblock In P.~Van~Dooren, S.~P. Bhattacharyya, R.~H. Chan, V.~Olshevsky, and
  A.~Routray, editors, {\em Numerical Linear Algebra in Signals, Systems and
  Control}, volume~80 of {\em Lect. Notes Electr. Eng.}, pages 55--76.
  Springer, Dodrecht, 2011.
\newblock \href {https://doi.org/10.1007/978-94-007-0602-6_3}
  {\path{doi:10.1007/978-94-007-0602-6_3}}.

\bibitem{BenB16}
P.~Benner and Z.~Bujanovi{\'c}.
\newblock On the solution of large-scale algebraic {R}iccati equations by using
  low-dimensional invariant subspaces.
\newblock {\em Linear Algebra Appl.}, 488:430--459, 2016.
\newblock \href {https://doi.org/10.1016/j.laa.2015.09.027}
  {\path{doi:10.1016/j.laa.2015.09.027}}.

\bibitem{BenBKetal18}
P.~Benner, Z.~Bujanovi{\'c}, P.~K{\"u}rschner, and J.~Saak.
\newblock {RADI}: a low-rank {ADI}-type algorithm for large scale algebraic
  {R}iccati equations.
\newblock {\em Numer. Math.}, 138(2):301--330, 2018.
\newblock \href {https://doi.org/10.1007/s00211-017-0907-5}
  {\path{doi:10.1007/s00211-017-0907-5}}.

\bibitem{BenBKetal20}
P.~Benner, Z.~Bujanovi{\'c}, P.~K{\"u}rschner, and J.~Saak.
\newblock A numerical comparison of different solvers for large-scale,
  continuous-time algebraic {R}iccati equations and {LQR} problems.
\newblock {\em {SIAM} J. Sci. Comput.}, 42(2):A957--A996, 2020.
\newblock \href {https://doi.org/10.1137/18M1220960}
  {\path{doi:10.1137/18M1220960}}.

\bibitem{BenEQetal14}
P.~Benner, P.~Ezzatti, Quintana-Ort{\'\i}~E. S., and A.~Rem{\'o}n.
\newblock A factored variant of the {N}ewton iteration for the solution of
  algebraic {R}iccati equations via the matrix sign function.
\newblock {\em Numer. Algorithms}, 66(2):363--377, 2014.
\newblock \href {https://doi.org/10.1007/s11075-013-9739-2}
  {\path{doi:10.1007/s11075-013-9739-2}}.

\bibitem{BenH17b}
P.~Benner and J.~Heiland.
\newblock Nonlinear feedback stabilization of incompressible flows via updated
  {R}iccati-based gains.
\newblock In {\em 2017 IEEE 56th Annual Conference on Decision and Control
  (CDC)}, pages 1163--1168, 2017.
\newblock \href {https://doi.org/10.1109/CDC.2017.8263813}
  {\path{doi:10.1109/CDC.2017.8263813}}.

\bibitem{BenH20}
P.~Benner and J.~Heiland.
\newblock Equivalence of {R}iccati-based robust controller design for index-1
  descriptor systems and standard plants with feedthrough.
\newblock In {\em 2020 European Control Conference (ECC)}, pages 402--407,
  2020.
\newblock \href {https://doi.org/10.23919/ECC51009.2020.9143771}
  {\path{doi:10.23919/ECC51009.2020.9143771}}.

\bibitem{morBenHW21}
P.~Benner, J.~Heiland, and S.~W.~R. Werner.
\newblock Robust output-feedback stabilization for incompressible flows using
  low-dimensional {$\mathcal{H}_{\infty}$}-controllers.
\newblock e-print 2103.01608, arXiv, 2021.
\newblock math.OC.
\newblock URL: \url{https://arxiv.org/abs/2103.01608}.

\bibitem{morBenKS21}
P.~Benner, M.~K{\"o}hler, and J.~Saak.
\newblock Matrix equations, sparse solvers: {M-M.E.S.S.}-2.0.1---{P}hilosophy,
  features and application for (parametric) model order reduction.
\newblock In P.~Benner, T.~Breiten, H.~Fa{\ss}bender, M.~Hinze, T.~Stykel, and
  R.~Zimmermann, editors, {\em Model Reduction of Complex Dynamical Systems},
  volume 171 of {\em International Series of Numerical Mathematics}, pages
  369--392. Birkh{\"a}user, Cham, 2021.
\newblock \href {https://doi.org/10.1007/978-3-030-72983-7_18}
  {\path{doi:10.1007/978-3-030-72983-7_18}}.

\bibitem{BenKS13a}
P.~Benner, P.~K{\"u}rschner, and J.~Saak.
\newblock A reformulated low-rank {ADI} iteration with explicit residual
  factors.
\newblock {\em Proc. Appl. Math. Mech.}, 13(1):585--586, 2013.
\newblock \href {https://doi.org/10.1002/pamm.201310273}
  {\path{doi:10.1002/pamm.201310273}}.

\bibitem{BenKS14b}
P.~Benner, P.~K{\"u}rschner, and J.~Saak.
\newblock Self-generating and efficient shift parameters in {ADI} methods for
  large {L}yapunov and {S}ylvester equations.
\newblock {\em Electron. Trans. Numer. Anal.}, 43:142--162, 2014.
\newblock URL:
  \url{https://etna.mcs.kent.edu/volumes/2011-2020/vol43/abstract.php?vol=43&pages=142-162}.

\bibitem{BenLP08}
P.~Benner, J.-R. Li, and T.~Penzl.
\newblock Numerical solution of large-scale {L}yapunov equations, {R}iccati
  equations, and linear-quadratic optimal control problems.
\newblock {\em Numer. Lin. Alg. Appl.}, 15(9):755--777, 2008.
\newblock \href {https://doi.org/10.1002/nla.622} {\path{doi:10.1002/nla.622}}.

\bibitem{BenS13}
P.~Benner and J.~Saak.
\newblock Numerical solution of large and sparse continuous time algebraic
  matrix {R}iccati and {L}yapunov equations: a state of the art survey.
\newblock {\em GAMM-Mitt.}, 36(1):32--52, 2013.
\newblock \href {https://doi.org/10.1002/gamm.201310003}
  {\path{doi:10.1002/gamm.201310003}}.

\bibitem{BenS14}
P.~Benner and T.~Stykel.
\newblock Numerical solution of projected algebraic {R}iccati equations.
\newblock {\em {SIAM} J. Numer. Anal.}, 52(2):581--600, 2014.
\newblock \href {https://doi.org/10.1137/130923993}
  {\path{doi:10.1137/130923993}}.

\bibitem{morBenS17}
P.~Benner and T.~Stykel.
\newblock Model order reduction for differential-algebraic equations: {A}
  survey.
\newblock In A.~Ilchmann and T.~Reis, editors, {\em Surveys in
  Differential-Algebraic Equations IV}, Differential-Algebraic Equations Forum,
  pages 107--160. Springer, Cham, 2017.
\newblock \href {https://doi.org/10.1007/978-3-319-46618-7_3}
  {\path{doi:10.1007/978-3-319-46618-7_3}}.

\bibitem{morBenW18}
P.~Benner and S.~W.~R. Werner.
\newblock Model reduction of descriptor systems with the {MORLAB} toolbox.
\newblock {\em IFAC-PapersOnLine}, 51(2):547--552, 2018.
\newblock \href {https://doi.org/10.1016/j.ifacol.2018.03.092}
  {\path{doi:10.1016/j.ifacol.2018.03.092}}.

\bibitem{morBenW19b}
P.~Benner and S.~W.~R. Werner.
\newblock {MORLAB} -- {Model Order Reduction LABoratory} (version 5.0), August
  2019.
\newblock See also: \url{https://www.mpi-magdeburg.mpg.de/projects/morlab}.
\newblock \href {https://doi.org/10.5281/zenodo.3332716}
  {\path{doi:10.5281/zenodo.3332716}}.

\bibitem{morBenW21c}
P.~Benner and S.~W.~R. Werner.
\newblock {MORLAB}---{T}he {M}odel {O}rder {R}eduction {LAB}oratory.
\newblock In P.~Benner, T.~Breiten, H.~Fa{\ss}bender, M.~Hinze, T.~Stykel, and
  R.~Zimmermann, editors, {\em Model Reduction of Complex Dynamical Systems},
  volume 171 of {\em International Series of Numerical Mathematics}, pages
  393--415. Birkh{\"a}user, Cham, 2021.
\newblock \href {https://doi.org/10.1007/978-3-030-72983-7_19}
  {\path{doi:10.1007/978-3-030-72983-7_19}}.

\bibitem{Del07}
M.~C. Delfour.
\newblock Linear quadratic differential games: {S}addle point and {R}iccati
  differential equation.
\newblock {\em {SIAM} J. Control Optim.}, 46(2):750--774, 2007.
\newblock \href {https://doi.org/10.1137/050639089}
  {\path{doi:10.1137/050639089}}.

\bibitem{morDesP82}
U.~B. Desai and D.~Pal.
\newblock A realization approach to stochastic model reduction and balanced
  stochastic realizations.
\newblock In {\em 21st IEEE Conference on Decision and Control}, pages
  1105--1112, 1982.
\newblock \href {https://doi.org/10.1109/CDC.1982.268322}
  {\path{doi:10.1109/CDC.1982.268322}}.

\bibitem{DoyGKF89}
J.~Doyle, K.~Glover, P.~P. Khargonekar, and B.~A. Francis.
\newblock State-space solutions to standard $\mathcal{H}_2$ and
  $\mathcal{H}_{\infty}$ control problems.
\newblock {\em {IEEE} Trans. Autom. Control}, 34(8):831--847, 1989.
\newblock \href {https://doi.org/10.1109/9.29425} {\path{doi:10.1109/9.29425}}.

\bibitem{morFreRM08}
F.~Freitas, J.~Rommes, and N.~Martins.
\newblock {G}ramian-based reduction method applied to large sparse power system
  descriptor models.
\newblock {\em {IEEE} Trans. Power Syst.}, 23(3):1258--1270, 2008.
\newblock \href {https://doi.org/10.1109/TPWRS.2008.926693}
  {\path{doi:10.1109/TPWRS.2008.926693}}.

\bibitem{GolV13}
G.~H. Golub and C.~F. Van~Loan.
\newblock {\em Matrix Computations}.
\newblock Johns Hopkins Studies in the Mathematical Sciences. Johns Hopkins
  University Press, Baltimore, fourth edition, 2013.

\bibitem{Hei16}
J.~Heiland.
\newblock A differential-algebraic {R}iccati equation for applications in flow
  control.
\newblock {\em {SIAM} J. Control Optim.}, 54(2):718--739, 2016.
\newblock \href {https://doi.org/10.1137/151004963}
  {\path{doi:10.1137/151004963}}.

\bibitem{supHeiW21a}
J.~Heiland and S.~W.~R. Werner.
\newblock Code, data and results for numerical experiments in ``{R}obust
  output-feedback stabilization for incompressible flows using low-dimensional
  $\mathcal{H}_{\infty}$-controllers'' (version 2.0), October 2021.
\newblock \href {https://doi.org/10.5281/zenodo.5532539}
  {\path{doi:10.5281/zenodo.5532539}}.

\bibitem{morHeiSS08}
M.~Heinkenschloss, D.~C. Sorensen, and K.~Sun.
\newblock Balanced truncation model reduction for a class of descriptor systems
  with application to the {O}seen equations.
\newblock {\em {SIAM} J. Sci. Comput.}, 30(2):1038--1063, 2008.
\newblock \href {https://doi.org/10.1137/070681910}
  {\path{doi:10.1137/070681910}}.

\bibitem{HeyJ09}
M.~Heyouni and K.~Jbilou.
\newblock An extended block {A}rnoldi algorithm for large-scale solutions of
  the continuous-time algebraic {R}iccati equation.
\newblock {\em Electron. Trans. Numer. Anal.}, 33:53--62, 2009.
\newblock URL:
  \url{https://etna.math.kent.edu/volumes/2001-2010/vol33/abstract.php?vol=33&pages=53-62}.

\bibitem{morJonS83}
E.~A. Jonckheere and L.~M. Silverman.
\newblock A new set of invariants for linear systems--application to reduced
  order compensator design.
\newblock {\em {IEEE} Trans. Autom. Control}, 28(10):953--964, 1983.
\newblock \href {https://doi.org/10.1109/TAC.1983.1103159}
  {\path{doi:10.1109/TAC.1983.1103159}}.

\bibitem{Kle68}
D.~L. Kleinman.
\newblock On an iterative technique for {R}iccati equation computations.
\newblock {\em {IEEE} Trans. Autom. Control}, 13(1):114--115, 1968.
\newblock \href {https://doi.org/10.1109/TAC.1968.1098829}
  {\path{doi:10.1109/TAC.1968.1098829}}.

\bibitem{Kue16}
P.~K{\"u}rschner.
\newblock {\em Efficient Low-Rank Solution of Large-Scale Matrix Equations}.
\newblock {Dissertation}, Department of Mathematics, Otto von Guericke
  University, Magdeburg, Germany, 2016.
\newblock URL: \url{http://hdl.handle.net/11858/00-001M-0000-0029-CE18-2}.

\bibitem{LanR95}
P.~Lancaster and L.~Rodman.
\newblock {\em Algebraic {R}iccati Equations}.
\newblock Oxford Science Publications. The Clarendon Press, Oxford University
  Press, New York, 1995.

\bibitem{LanFA07}
A.~Lanzon, Y.~Feng, and B.~D.~O. Anderson.
\newblock An iterative algorithm to solve algebraic {R}iccati equations with an
  indefinite quadratic term.
\newblock In {\em 2007 European Control Conference (ECC)}, pages 3033--3039,
  2007.
\newblock \href {https://doi.org/10.23919/ecc.2007.7068239}
  {\path{doi:10.23919/ecc.2007.7068239}}.

\bibitem{LanFAetal08}
A.~Lanzon, Y.~Feng, B.~D.~.O. Anderson, and M.~Rotkowitz.
\newblock Computing the positive stabilizing solution to algebraic {R}iccati
  equations with an indefinite quadratic term via a recursive method.
\newblock {\em {IEEE} Trans. Autom. Control}, 53(10):2280--2291, 2008.
\newblock \href {https://doi.org/10.1109/TAC.2008.2006108}
  {\path{doi:10.1109/TAC.2008.2006108}}.

\bibitem{Lau79}
A.~J. Laub.
\newblock A {S}chur method for solving algebraic {R}iccati equations.
\newblock {\em {IEEE} Trans. Autom. Control}, 24(6):913--921, 1979.
\newblock \href {https://doi.org/10.1109/TAC.1979.1102178}
  {\path{doi:10.1109/TAC.1979.1102178}}.

\bibitem{Lei04}
F.~Leibfritz.
\newblock {$COMPl_{e}ib$}: \emph{CO}nstrained \emph{M}atrix-optimization
  \emph{P}roblem \emph{lib}rary -- a collection of test examples for nonlinear
  semidefinite programs, control system design and related problems.
\newblock Tech.-report, University of Trier, 2004.
\newblock URL:
  \url{http://www.friedemann-leibfritz.de/COMPlib_Data/COMPlib_Main_Paper.pdf}.

\bibitem{LiW02}
J.-R. Li and J.~White.
\newblock Low rank solution of {L}yapunov equations.
\newblock {\em {SIAM} J. Matrix Anal. Appl.}, 24(1):260--280, 2002.
\newblock \href {https://doi.org/10.1137/S0895479801384937}
  {\path{doi:10.1137/S0895479801384937}}.

\bibitem{LinS15}
Y.~Lin and V.~Simoncini.
\newblock A new subspace iteration method for the algebraic {R}iccati equation.
\newblock {\em Numer. Linear Algebra Appl.}, 22(1):26--47, 2015.
\newblock \href {https://doi.org/10.1002/nla.1936}
  {\path{doi:10.1002/nla.1936}}.

\bibitem{Loc01}
A.~Locatelli.
\newblock {\em Optimal Control: {A}n Introduction}.
\newblock Birkh{\"a}user, Basel, 2001.

\bibitem{McFG90}
D.~C. McFarlane and K.~Glover.
\newblock {\em Robust Controller Design Using Normalized Coprime Factor Plant
  Descriptions}, volume 138 of {\em Lect. Notes Control Inf. Sci.}
\newblock Springer, Berlin, Heidelberg, 1990.
\newblock \href {https://doi.org/10.1007/BFB0043199}
  {\path{doi:10.1007/BFB0043199}}.

\bibitem{morMoeRS11}
J.~M{\"o}ckel, T.~Reis, and T.~Stykel.
\newblock Linear-quadratic {G}aussian balancing for model reduction of
  differential-algebraic systems.
\newblock {\em Internat. J. Control}, 84(10):1627--1643, 2011.
\newblock \href {https://doi.org/10.1080/00207179.2011.622791}
  {\path{doi:10.1080/00207179.2011.622791}}.

\bibitem{morMusG91}
D.~Mustafa and K.~Glover.
\newblock Controller reduction by $\mathcal{H}_\infty$-balanced truncation.
\newblock {\em {IEEE} Trans. Autom. Control}, 36(6):668--682, 1991.
\newblock \href {https://doi.org/10.1109/9.86941} {\path{doi:10.1109/9.86941}}.

\bibitem{morOpdJ88}
P.~C. Opdenacker and E.~A. Jonckheere.
\newblock A contraction mapping preserving balanced reduction scheme and its
  infinity norm error bounds.
\newblock {\em {IEEE} Trans. Circuits Syst.}, 35(2):184--189, 1988.
\newblock \href {https://doi.org/10.1109/31.1720} {\path{doi:10.1109/31.1720}}.

\bibitem{morRob80}
J.~D. Roberts.
\newblock Linear model reduction and solution of the algebraic {R}iccati
  equation by use of the sign function.
\newblock {\em Internat. J. Control}, 32(4):677--687, 1980.
\newblock Reprint of Technical Report No. TR-13, CUED/B-Control, Cambridge
  University, Engineering Department, 1971.
\newblock \href {https://doi.org/10.1080/00207178008922881}
  {\path{doi:10.1080/00207178008922881}}.

\bibitem{Saa09}
J.~Saak.
\newblock {\em Efficient Numerical Solution of Large Scale Algebraic Matrix
  Equations in {PDE} Control and Model Order Reduction}.
\newblock {Dissertation}, Department of Mathematics, University of Technology
  Chemnitz, Chemnitz, Germany, 2009.
\newblock URL: \url{https://nbn-resolving.org/urn:nbn:de:bsz:ch1-200901642}.

\bibitem{SaaKB21-mmess-2.1}
J.~Saak, M.~K{\"o}hler, and P.~Benner.
\newblock {M-M.E.S.S.} -- {T}he {M}atrix {E}quations {S}parse {S}olvers library
  (version 2.1), April 2021.
\newblock See also:~\url{https://www.mpi-magdeburg.mpg.de/projects/mess}.
\newblock \href {https://doi.org/10.5281/zenodo.4719688}
  {\path{doi:10.5281/zenodo.4719688}}.

\bibitem{morSaaV18}
J.~Saak and M.~Voigt.
\newblock Model reduction of constrained mechanical systems in {M-M.E.S.S.}
\newblock {\em IFAC-PapersOnLine}, 51(2):661--666, 2018.
\newblock \href {https://doi.org/10.1016/j.ifacol.2018.03.112}
  {\path{doi:10.1016/j.ifacol.2018.03.112}}.

\bibitem{San74}
N.~Sandell.
\newblock On {N}ewton's method for {R}iccati equation solution.
\newblock {\em {IEEE} Trans. Autom. Control}, 19(3):254--255, 1974.
\newblock \href {https://doi.org/10.1109/TAC.1974.1100536}
  {\path{doi:10.1109/TAC.1974.1100536}}.

\bibitem{Sim16}
V.~Simoncini.
\newblock Analysis of the rational {K}rylov subspace projection method for
  large-scale algebraic {R}iccati equations.
\newblock {\em {SIAM} J. Matrix Anal. Appl.}, 37(4):1655--1674, 2016.
\newblock \href {https://doi.org/10.1137/16M1059382}
  {\path{doi:10.1137/16M1059382}}.

\bibitem{Son98}
E.~D. Sontag.
\newblock {\em Mathematical Control Theory}, volume~6 of {\em Texts in Applied
  Mathematics}.
\newblock Springer, New York, second edition, 1998.

\bibitem{Sti18}
T.~Stillfjord.
\newblock Singular value decay of operator-valued differential {L}yapunov and
  {R}iccati equations.
\newblock {\em SIAM J. Control Optim.}, 56(5):3598--3618, 2018.
\newblock \href {https://doi.org/10.1137/18M1178815}
  {\path{doi:10.1137/18M1178815}}.

\bibitem{Sty02}
T.~Stykel.
\newblock {\em Analysis and Numerical Solution of Generalized {L}yapunov
  Equations}.
\newblock {Dissertation}, Faculty II -- Mathematics and Nature Sciences,
  Technische Universit{\"a}t Berlin, Berlin, Germany, 2002.
\newblock \href {https://doi.org/10.14279/depositonce-578}
  {\path{doi:10.14279/depositonce-578}}.

\bibitem{Sty08}
T.~Stykel.
\newblock Low-rank iterative methods for projected generalized {L}yapunov
  equations.
\newblock {\em Electron. Trans. Numer. Anal.}, 30:187--202, 2008.
\newblock URL:
  \url{https://etna.math.kent.edu/volumes/2001-2010/vol30/abstract.php?vol=30&pages=187-202}.

\bibitem{Var95b}
A.~Varga.
\newblock On computing high accuracy solutions of a class of {R}iccati
  equations.
\newblock {\em Control--Theory and Advanced Technology}, 10(4):2005--2016,
  1995.

\bibitem{Wei16}
H.~K. Weichelt.
\newblock {\em Numerical Aspects of Flow Stabilization by {R}iccati Feedback}.
\newblock {D}issertation, Department of Mathematics, Otto von Guericke
  University, Magdeburg, Germany, 2016.
\newblock \href {https://doi.org/10.25673/4493} {\path{doi:10.25673/4493}}.

\bibitem{ZhaFC20}
L.~Zhang, H.-Y. Fan, and E.~K. Chu.
\newblock Inheritance properties of {K}rylov subspace methods for
  continuous-time algebraic {R}iccati equations.
\newblock {\em J. Comput. Appl. Math.}, 371:112685, 2020.
\newblock \href {https://doi.org/10.1016/j.cam.2019.112685}
  {\path{doi:10.1016/j.cam.2019.112685}}.

\end{thebibliography}

\end{document}